\documentclass[11pt]{article}
\usepackage{graphicx}
\usepackage{amssymb,amsmath,amsthm}
\usepackage{color}
\usepackage{diagrams}
\newtheorem{theorem}{Theorem}[section]
\newtheorem{corollary}[theorem]{Corollary}
\newtheorem{lemma}[theorem]{Lemma}
\theoremstyle{definition}
\newtheorem{definition}[theorem]{Definition}
\theoremstyle{remark}
\newtheorem{remark}[theorem]{Remark}
\theoremstyle{definition}
\newtheorem{example}[theorem]{Example}

\begin{document}
\title{A New Perspective for an Existing Homology Theory of Links Embedded in I-Bundles}
\author{Jeffrey Boerner\footnote{partially supported by the University of Iowa Department of Mathematics NSF VIGRE grant DMS-0602242}\\
University of Iowa}

\renewcommand{\today}{}
\maketitle

\begin{abstract}  This paper introduces a homology theory for links in I-bundles over an orientable surface.  The theory is unique in that the elements of the chain groups are surfaces instead of diagrams.  It is then shown this theory yields the same results as the homology theory constructed by Asaeda, Przytycki and Sikora.

\end{abstract}

\section {Introduction}

%Historical Perspective

In [K] Khovanov introduced a homology theory for links in $S^3$ that was a categorification of the Jones polynomial.  In [APS] Asaeda, Przytycki and Sikora extended this theory to links embedded in I-bundles.  Their homology theory incorporated some of the topology of the I-bundle into their invariant.  

%How this relates to the other work

Turner and Turaev showed in [T] that the homology from [APS] could be recreated using embedded surfaces as elements of the chain groups instead of decorated diagrams.  In this paper we accomplish that in the case of I-bundles over orientable surfaces.

%Structure of the paper

Section 2 contains definitions and explains the skein relations on surfaces that are used.  Section 3 defines the grading on the chain groups and which surfaces generate the chain groups.

The boundary operator is defined in section 4 and it is also shown that it is well-defined with respect to the relations.  In section 5 it is proved that the boundary operator squared is equal to zero, and thus the boundary operator together with the chain groups form a chain complex.

Finally, in section 6 it is shown that the homology produced from the chain complex coincides with the homology from [APS]

\section {Definitions}

\begin{definition}

Let $S$ be a surface properly embedded in a 3-manifold $N$.  A boundary circle of $S$ is said to be \textbf{inessential} if it bounds a disk in $N$, otherwise it is said to be \textbf{essential}.

\end{definition}

\begin{definition}

If $S$ is an oriented surface and $c$ is an oriented boundary component of $S$ then the orientation of $S$ is \textbf{compatible} with the orientation of $c$ if the boundary orientation of $c$ from $S$ agrees with the orientation of $c$.    Two oriented boundary curves of an orientable connected surface are \textbf{compatible} if both curves are compatible with the same orientation on the surface.

\end{definition}

\begin{remark}
If $S$ is a connected unoriented orientable surface and $c$ is an oriented boundary component of $S$ then there is exactly one orientation for the other boundary curves to be oriented compatibly with $c$.
\end{remark}

\begin{center}

\begin{tabular}{ccccc}
\includegraphics[height = .8 in]{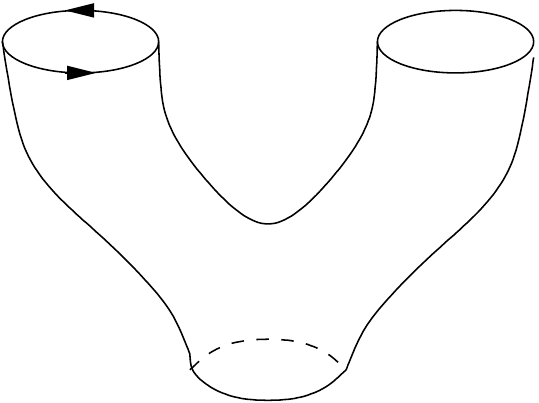}  & \includegraphics[height = .8 in]{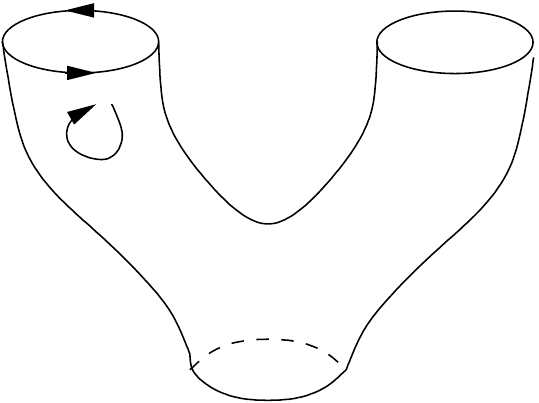}     & \includegraphics[height = .8 in]{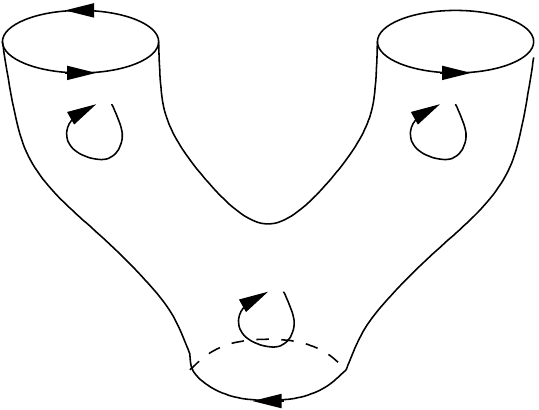}  \\

Boundary is oriented    &   Compatible orientation on surface     &  Rest of boundary is oriented  \\
&& compatibly
\end{tabular}

\end{center}

%\end{definition}

Let $F$ be an orientable compact surface not necessarily with boundary.  Let $D$ be a link diagram in $F$, such that the crossings of $D$ are enumerated.

\begin{definition}

A \textbf{state} of the diagram $D$ is a choice of smoothing at each crossing.  Thus a state is represented by a collection of disjoint simple closed curves in the surface $F$.

\end{definition}

\begin{definition}

A \textbf{state surface} with respect to the diagram $D$ has the following properties:

\begin{itemize}

\item A state surface is an orientable compact surface properly embedded in $F$ x $I$.
\item A state surface has a state of $D$ as its boundary in the top ($F$ x \{0\}) and essential oriented circles as its boundary in the bottom ($F$ x \{1\}).  
\item Inessential boundary curves of state surfaces are not oriented, but essential boundary curves in the top may or may not be oriented.  
\item If one component of a state surface has an oriented essential boundary curve, then all essential boundary curves on that component must be oriented compatibly.   
\item State surfaces may be marked with dots.
\end{itemize}

Two state surfaces are equivalent if they are isotopic relative to the boundary.  Thus the dots on the state surfaces are allowed to move freely within components but dots may not switch components.  

\end{definition}

Let $M$ be the free $\mathbb{Z}$-module generated by state surfaces with respect to the diagram $D$.

\medskip

\begin{remark}

In order to continue, local relations need to be defined on $M$.  The manner that this is done is to define a submodule of $M$, $B$, in order that the relations hold in $M/B$.  Thus if we want $C=D$, then we have $C-D$ as a generator of $B$, so then in $M/B$, $C$ is equivalent to $D$.  These are skein relations, so if a relation is $P=P'$ it means that $S=S'$ in $M/B$ if there is a 3-ball, $A$, in the manifold such that $S$ and $S'$ agree outside of $A$ while $S\cap A=P$ and $S'\cap A=P'$.

Also if $z\in \mathbb{Z}$, then the relation $P=z$ means if a surface has $P$ as a subsurface, then the original surface is equal to $z$ times the surface where $P$ is removed in $M/B$.

\begin{example}

Note the Neck-Cutting Relation
\begin{center}
\includegraphics[height=1 in]{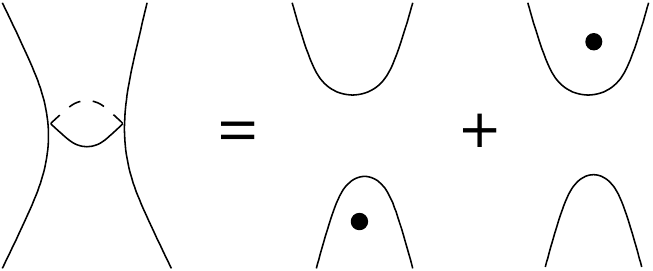}
\end{center}

implies the equality

\begin{center}
\includegraphics[height = 1.5 in]{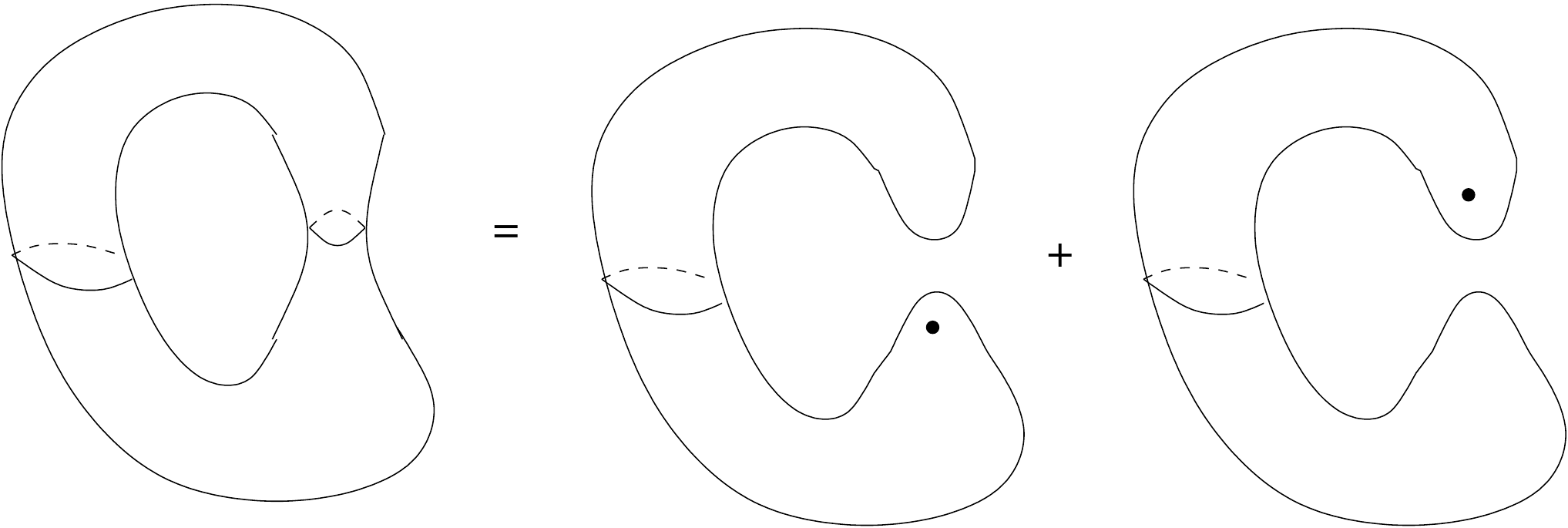}
\end{center}

and the relation 

\begin{center}
\includegraphics[height = 1 in]{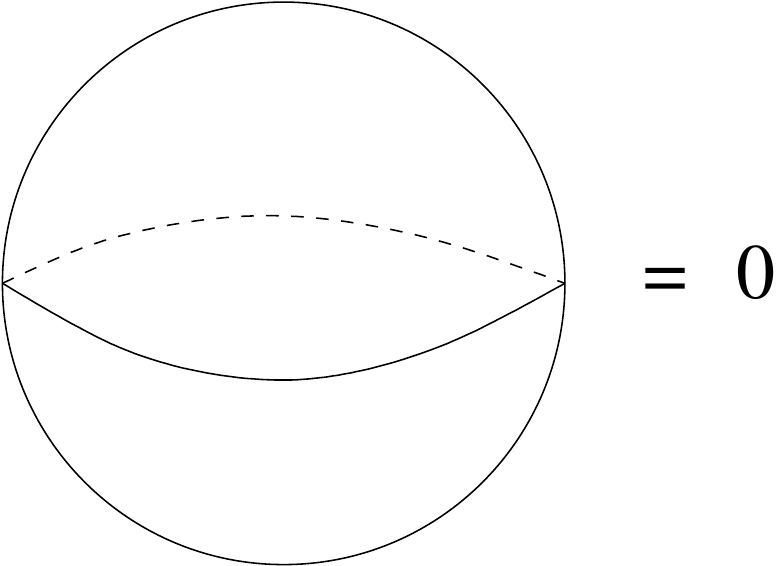}
\end{center}

implies 

\begin{center}
\includegraphics[height = 1 in]{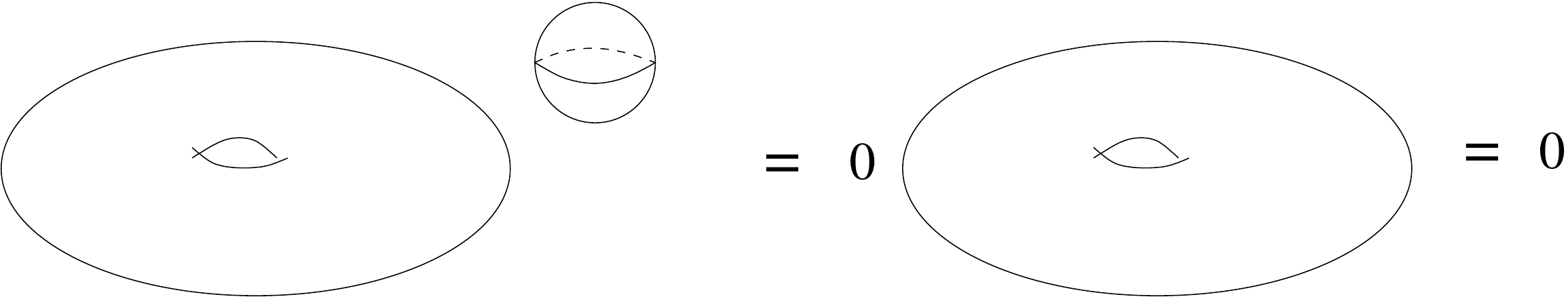}
\end{center}

\end{example}

\end{remark}

\begin{definition} Let \textbf{B} be the submodule generated by the following relations: 
\end{definition}

\begin{itemize}

\item Neck-Cutting Relation (NC)
\begin{center}
\includegraphics[height=1 in]{BNrel.pdf}
\end{center}

\item Sphere bounding a ball equals zero 
\begin{center}
\includegraphics[height=1 in]{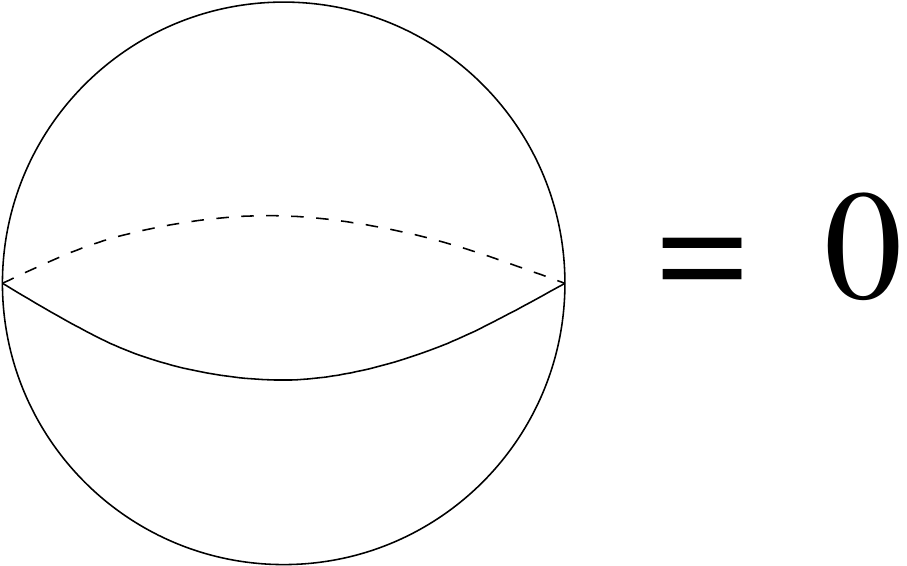}
\end{center}

\item Sphere with a dot bounding a ball equals one (SD)
\begin{center}
\includegraphics[height=1 in]{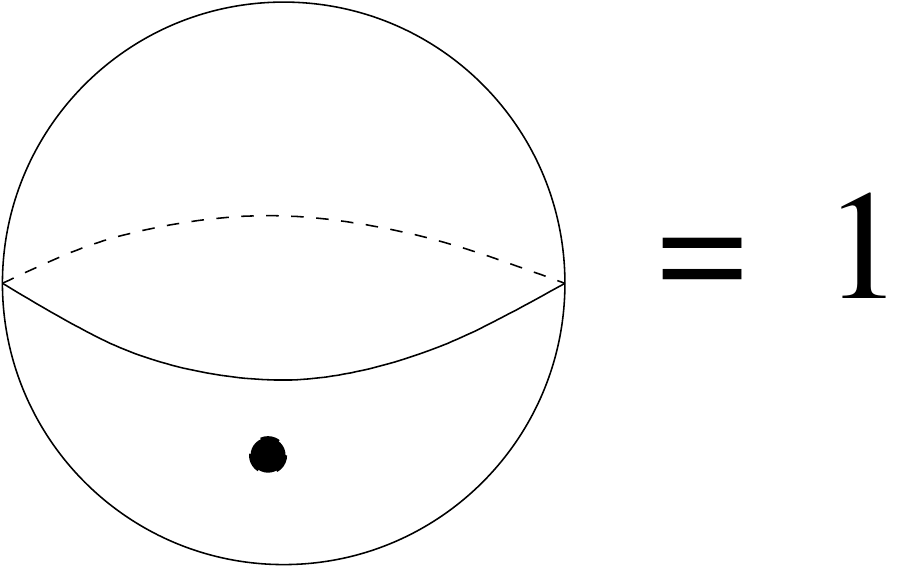}
\end{center}

\item A component with two dots equals zero.
\begin{center}
\includegraphics[height=.6 in]{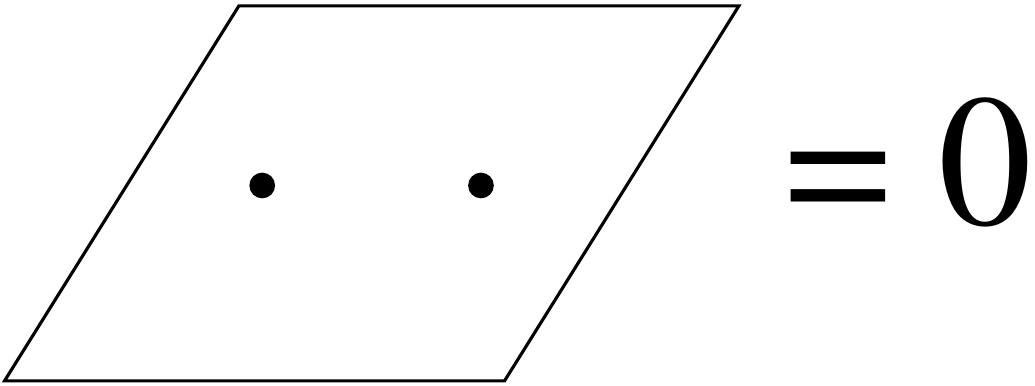}
\end{center}

\item An annulus with its boundary completely in the bottom equals one. (BDA)

\item A surface with a non-disk, non-sphere component, that has a dot on that component equals zero.    (NDD)

\item A surface with an incompressible component with negative euler characteristic equals zero.   (NEC)

\item If the annulus in the figure below on the left side of the inequality is incompressible then we have the relation:

\begin{center} \includegraphics[height=1 in]{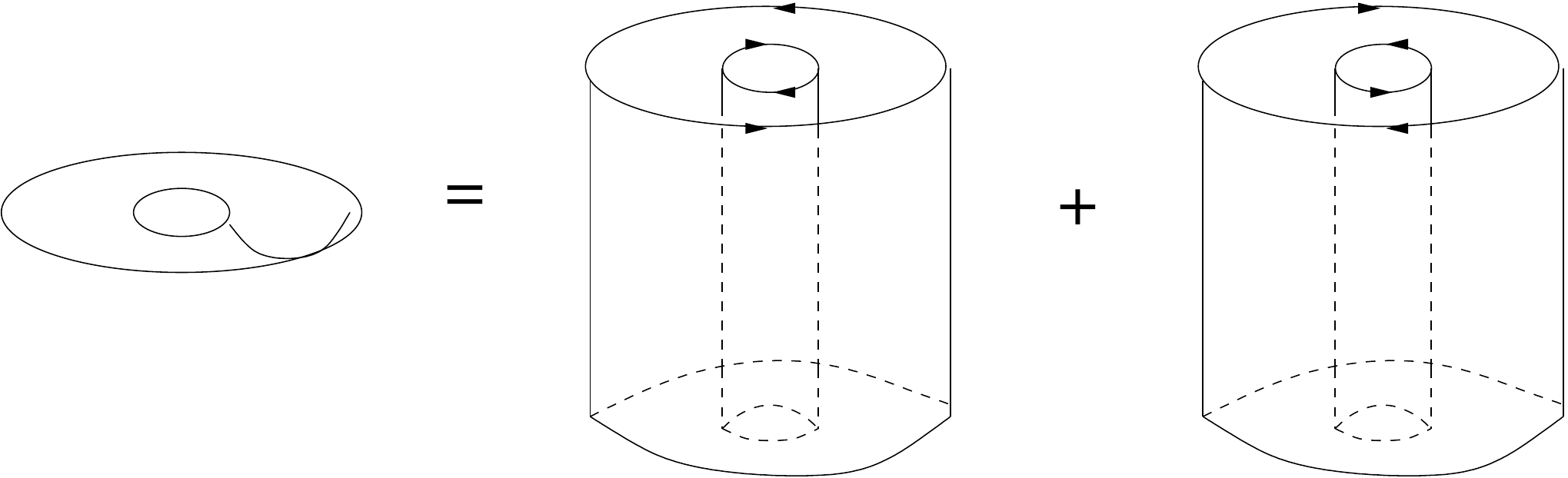}      \hspace{.05 in} \raisebox{.5 in}{(UTA)} \end{center}

\end{itemize}

Let $P= M/B$.  The elements of P will be referred to as foams.

\section{Chain Groups}

\begin{definition}

Let $S$ be a state surface.

$I(S) = \{\#$ of positive smoothings in the state corresponding to the top boundary of $S$\} $ -  $ \{\# of negative smoothings in the state corresponding to the top boundary of $S$\}

$J(S) = I(S) + 2(2d  - \chi(S))$  where $d$ is the number of dots on $S$.

\end{definition}

\medskip

The third index corresponds to the oriented essential disjoint simple closed curves in the bottom of $F$ x $I$.  Note that given two parallel oriented simple closed curves on a surface it is not difficult to determine if their orientations agree or disagree.  To define the third index it is necessary to specify that one of the orientations possible is the positive one and the opposite orientation is the negative one.  Thus for each homotopy class of simple closed curve we have chosen a positive orientation and a negative orientation.

\begin{definition}
Let $\gamma_1,\dots,\gamma_n$ be a family of disjoint simple closed curves in the bottom of a state surface $S$.    If $\gamma_i$ and $\gamma_j$ are parallel then $\gamma_i=\gamma_j$.
Then 

$K(S) = \sum_{i=1}^n k_i\gamma_i$, where 

$$   
k_i  = \left\{ 
\begin{array}{ccc}
1,&\mbox{ if } \text{$\gamma_i$ is oriented in the positive direction}                  \\
-1,&\mbox{ if } \text{$\gamma_i$ is oriented in the negative direction}                  \\
\end{array}\right.
$$

\end{definition}

\medskip

In order for these definitions to make sense they need to be well defined in the quotient.  Thus we have the following Lemma:

\begin{lemma}
The indices $I(S)$, $J(S)$ and $K(S)$ are well defined for $S \in P$.
\end{lemma}

\begin{proof}
Showing the index $I(S)$ is well-defined is immediate because the top circles are not affected by the relations.

In order to show the index $J(S)$ is well defined we need to consider the Neck-Cutting Relation, the SD relation, the UTA relation and the BDA relation.

To consider the Neck-Cutting Relation, first note that when a neck is cut the euler characteristic goes up by two.  Also, each summand adds a dot, thus $2d  - \chi(S)$ remains the same.  

A sphere has euler characteristic two, and note if the sphere has a dot, then $2d-\chi(S) = 2 -2 =0$, so removing a sphere with a dot does not affect the $J(S)$ index.  

When considering the UTA relation, note that both sides of the equality have annuli without dots which do not contribute dots or euler characteristic, thus the $J(S)$ grading on each side of the equality is the same.  

We may also remove an annulus without a dot, since an annulus has euler characteristic zero, so then $2d-\chi(S) = 0$.  Thus the BDA relation is well defined.

The $K(S)$ index is only dependent on the curves in the bottom, so we only need to consider the relations that affect the curves in the bottom.  The only relations that do this are the BDA relation and the UTA relation.  

Note that a bottom annulus has two essential homotopic curves.  Since these curves are in the bottom they must be oriented and they must be oriented in a way that is compatible with each other.  Thus they must be oriented in the opposite direction of one other.  Therefore in the $K(S)$ index, these curves cancel out in the sum, so the annulus does not contribute to the $K(S)$ index, so we may remove the annulus without affecting the grading.

For the UTA relation the left side of the equation has no curves in the bottom and in each summand on the right side there are homotopic curves oriented in the opposite direction.  Since these curves have opposite orientation, they cancel out in the sum that determines $K(S)$, so neither side of the equality contributes to the $K(S)$ grading.

\end{proof}

\medskip

Let $C_{i,j,s}(D)$ be the submodule of $P$ generated by all foams $S$ of disks and incompressible oriented vertical annuli such that  $I(S)=i$, $J(S)=j$ and $K(S)=s$.

\section {The Boundary Operator}
In order to define the boundary operator the following Lemma is needed.

\begin{lemma}
\label{nonorsurf}
Placing a bridge doesn't affect the number of boundary curves if and only if a non-orientable surface is created.

\end{lemma}

\begin{proof}
This is an if and only if proof, so first left to right will be shown and then right to left.
($\rightarrow$) Assume placing a bridge turns one boundary curve into one boundary curve.  

The only way this can happen is if the curves are as they are in the diagram, since if the upper left connected to the upper right placing a bridge results in two curves, and if the upper left connects to the lower left then we are starting with two curves.  

\begin{center}
\includegraphics[height=1.5 in]{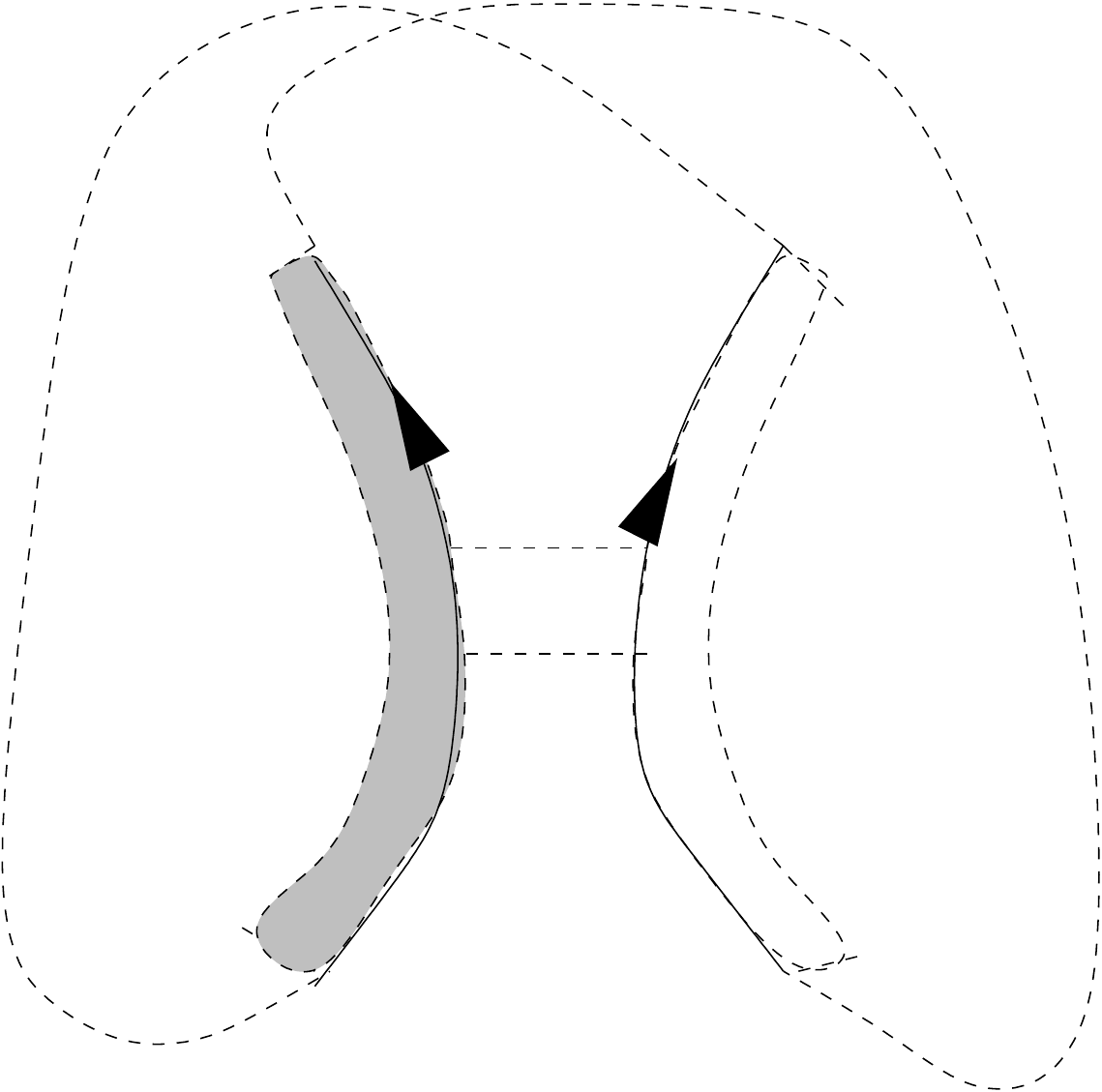}
\end{center}

Now note the boundary curve is a circle, so we can place an orientation on it from the boundary orientation of the surface.  Since the original surface is orientable we can color the side to the left of the circle (if we face the direction the arrow is pointing) and leave the other side blank.  Thus each side of the component is determined to be dark or blank by the orientation of the curve.

Now note when a bridge is placed it must connect a dark side to a blank side, thus resulting in a non-orientable surface.

($\leftarrow$) Assume a non-orientable surface is created as a result of placing a bridge.

We started with an orientable surface, thus before placing the bridge it was possible to color the surface with a dark side and a blank side.  If it were possible to place the bridge in such a way that the coloring remained consistent, then the resulting surface would be orientable.  Thus we can assume this is not possible and the bridge connects the dark side to the blank side.  Now consider the original surface before placing the bridge.  It is possible to orient the curve so that the dark side is on the left all the time since the original surface is orientable.  

\begin{center}
\includegraphics[height=1.5 in]{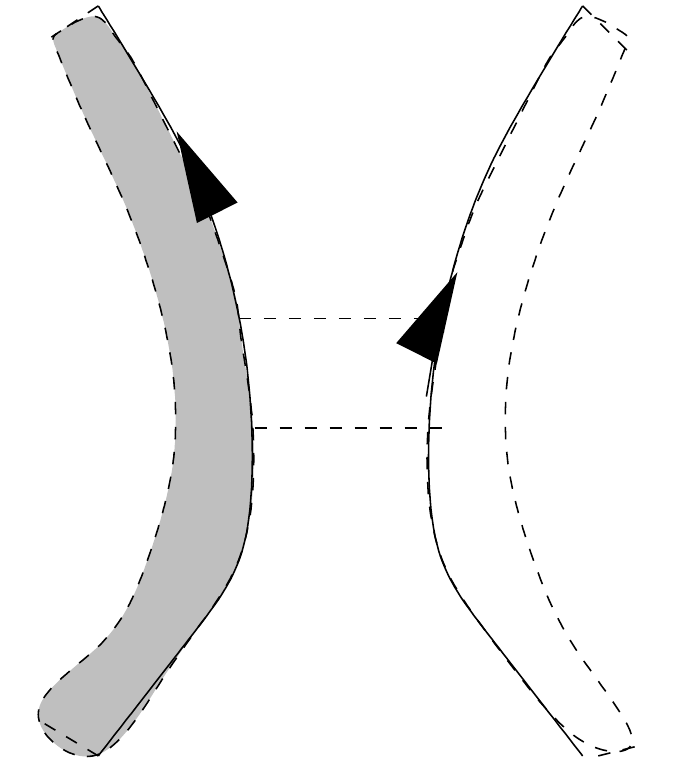} 
\end{center}

Then note for the orientation to be consistent on the circle, the top left connects to either the bottom left or the bottom right.  However, the bottom left results in two curves to start with, so the top left must connect to the bottom left.  Thus we are left with the following diagram and one can note the change of crossing results in one curve becoming one curve.

\begin{center}
\includegraphics[height=1.5 in]{nonorient.pdf}
\end{center}

\end{proof}

Let $p$ be a crossing of the diagram $D$.  We will define the partial boundary operator $d_p: C_{i,j,s}\rightarrow C_{i-2,j,s}$.  It is necessary to define $\bar{d}_p$:

\begin{center}
\includegraphics[height=1 in]{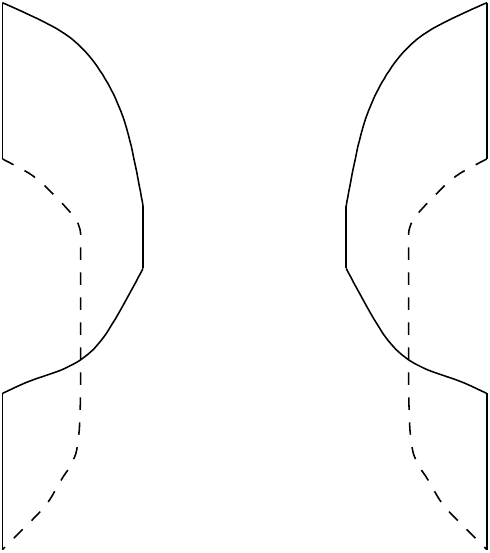}  \hspace{1 in}  \includegraphics[height=1 in]{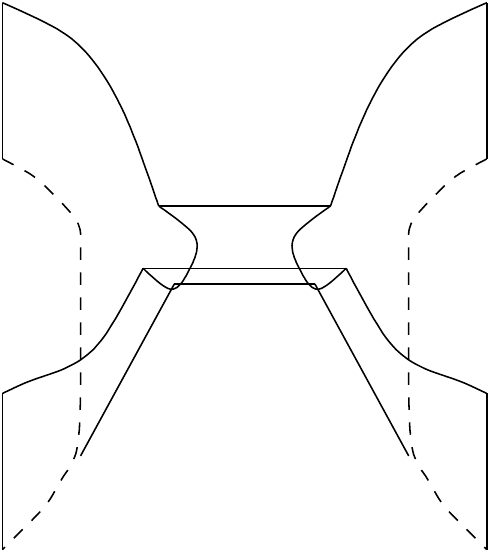}

$S$ at $p$-th crossing  \hspace{.7 in}   $\bar{d}_p(S)$ at $p$-th crossing
\end{center}

Informally this operation will be called placing a bridge at the $p$-th crossing.

\medskip

There are two situations when $\bar{d}_p$ cannot be applied: 

\begin{enumerate}

\item The orientations on boundary curves at $p$ are not compatible.  One example of this is in the figure below.

\begin{center}
\includegraphics[height = 1 in]{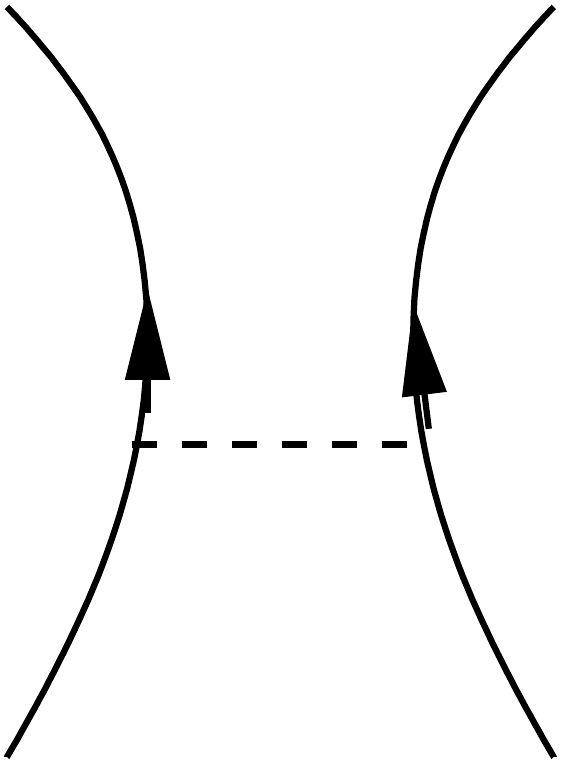}
\end{center}

Notationally this will be referred to as (EO).

\item Placing a bridge at $p$ creates a non-orientable surface.  This will be referred to as (NOS).

\end{enumerate}
 
$\bar{d}_p(S)$ is not defined if EO or NOS occurs at the $p$-th crossing of $S$.

\medskip

We must now address how to orient the boundary curves of $\bar{d}_p(S)$.  The general rule is to preserve the orientation of arcs that appear in both the boundary of $S$ and in the boundary of $\bar{d}_p(S)$ whenever possible, with the additional requirement that inessential circles of $\bar{d}_p(S)$ are not oriented.  This will be made more rigorous below. 

\medskip

Note that placing a bridge may either turn two boundary curves into one boundary curve, one boundary curve into two boundary curves, or one boundary curve into one boundary curve.  

\medskip

We will show how $\bar{d}_p$ behaves with respect to orientation by cases:

\begin{enumerate}

\item Two boundary curves become one boundary curve after placing a bridge.

If the original boundary curves are unoriented then the resulting boundary curve is unoriented as well.

If one of the original boundary curves is oriented the resulting boundary curve is oriented as in the figure below.

\begin{center}
\includegraphics[height = 1 in]{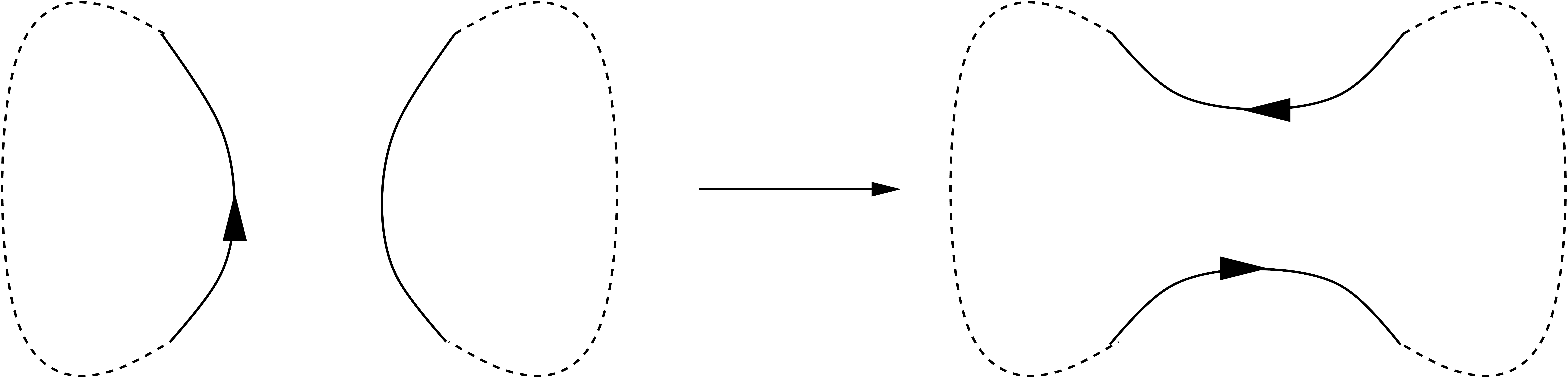}
\end{center}

Also, any essential curves that were unoriented are now oriented compatibly with the oriented curves on the same component.

If both of the original boundary curves are oriented then EO may occur.  If not then the orientations behave as in the figure.

\begin{center}
\includegraphics[height = 1 in]{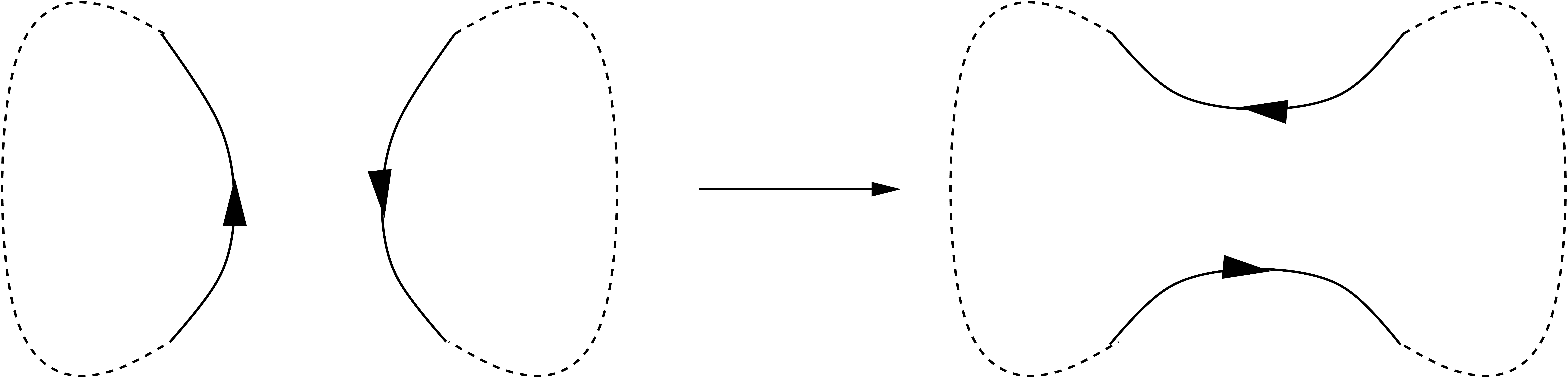}
\end{center}

\item One boundary curve becomes two boundary curves.

\begin{center}
\includegraphics[height = 1 in]{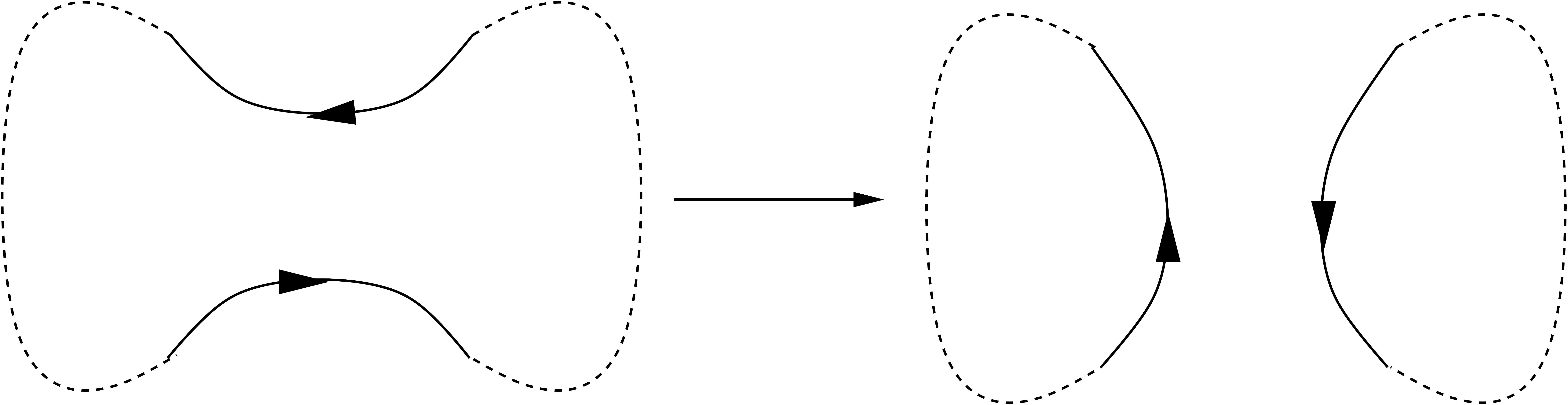}
\end{center}

\item One boundary curve becomes one boundary curve.

This only happens if NOS occurs by Lemma \ref{nonorsurf} and $\bar{d}_p(S)$ is not defined in that situation.

\end{enumerate}

Note inessential boundary curves of $\bar{d}_p(S)$ are not oriented.  Thus when evaluating $\bar{d}_p(S)$ by the rules above inessential boundary curves actually have no orientation in $\bar{d}_p(S)$.

\medskip

We define
$$
d_p(S) = \left\{ 
\begin{array}{ccc}
0,&\mbox{ if } \text{$p$-th crossing of $S$ is smoothed negatively}                  \\
0,&\mbox{ if } \text{(EO) occurs}                  \\
0, &\mbox{ if } \text{(NOS) occurs }                          \\
\text{$\bar{d}_p(S)$}, &\mbox{ else }
\end{array}\right.
$$

Then we have the differential $d: C_{i,j,s}(D) \rightarrow C_{i-2,j,s}(D)$, defined by 

\begin{equation}d(S) = 
\sum_{p \text{ a crossing of $D$}} (-1)^{t(S,p)}d_p(S),\end{equation}
where $t(S,p) = |$\{ $j$ a crossing of $D$ : $j$ is after $p$ in the ordering of crossings and $j$ is smoothed negatively in boundary state of S\}$|$

\medskip

\begin{lemma}
\label{vertann}
Let $S \in C_{i,j,s}(D)$ and p be a crossing of the diagram D.  Then $d_p(S)$ is not a vertical unoriented annulus or an oriented annulus with its boundary only in the top.
\end{lemma}

\begin{proof}
All boundary curves on the bottom are essential and oriented by the definition of a state surface.  A bridge does not affect the curves on the bottom thus they always stay essential and oriented, so a vertical unoriented annulus cannot occur.

We only need to consider a foam with components consisting of disks and oriented vertical annuli by the relations.  Note that placing a bridge between any of these components only creates connected components with one boundary curve in the top.  Thus we can conclude that to create annuli with their boundary only in the top, we need to bridge components to themselves.

Note an annulus could come from bridging a disk to itself, but a disk has no orientation on its boundary components.  Thus the resulting annulus would not have oriented boundary components.

Thus we can conclude that an oriented annulus with boundary in the top could only come from an oriented vertical annulus after applying the boundary operator.  Note, a vertical annulus has two boundary components.  Suppose after applying $d_p$ to this foam we now have two essential boundary components in the top, so three essential boundary components total.  The annulus we are attempting to create has both boundary components in the top, so in order to arrive at this annulus we would need to compress to get two components, one with two curves in the top and one with one curve in the bottom.  However, this is not possible, since the boundary curve is essential and it cannot be the only boundary curve on a surface.

\end{proof}

\begin{remark}

Here are two items to observe for the proof of the following Lemmas and Theorem:

\begin{enumerate}

\item Bridging an essential boundary curve to itself can produce at most one inessential boundary curve at a time.  This is due to the fact that if two inessential boundary curves are created from one curve, then the original curve was also inessential.

\item If a component has an inessential boundary component then either this component is a disk or it is compressible. (Just push the disk the curve bounds into $F$ x $I$ to obtain a compressing disk.)

\end{enumerate}

\end{remark}

\begin{lemma}
\label{nontrivdot}
A foam with a component that has an essential boundary component and also has a dot is trivial in the quotient.

\end{lemma}

\begin{proof}
If this component is incompressible we are done.  If not, then compress this component.  If the compressing disk was non-separating, then the result of compressing is a component with two dots, thus it is trivial in the quotient.

If the compressing disk is separating we end up with two components and one of them has an essential boundary curve, thus one of the components is not a disk, but they both have dots.

We are now left with two components, each with dots, and one of them has an essential boundary component. 

Note the component that has an essential boundary curve and a dot satisfies the assumptions of the Lemma.  Thus we may continue in this manner.  Note that each time the neck-cutting relation is applied the euler characteristic goes up by two and since we are dealing with compact surfaces the euler characteristic is bounded above.  Thus this process terminates and we eventually arrive at an incompressible component with a dot that is not a disk or a sphere.

\begin{center}
\includegraphics[height = 2 in]{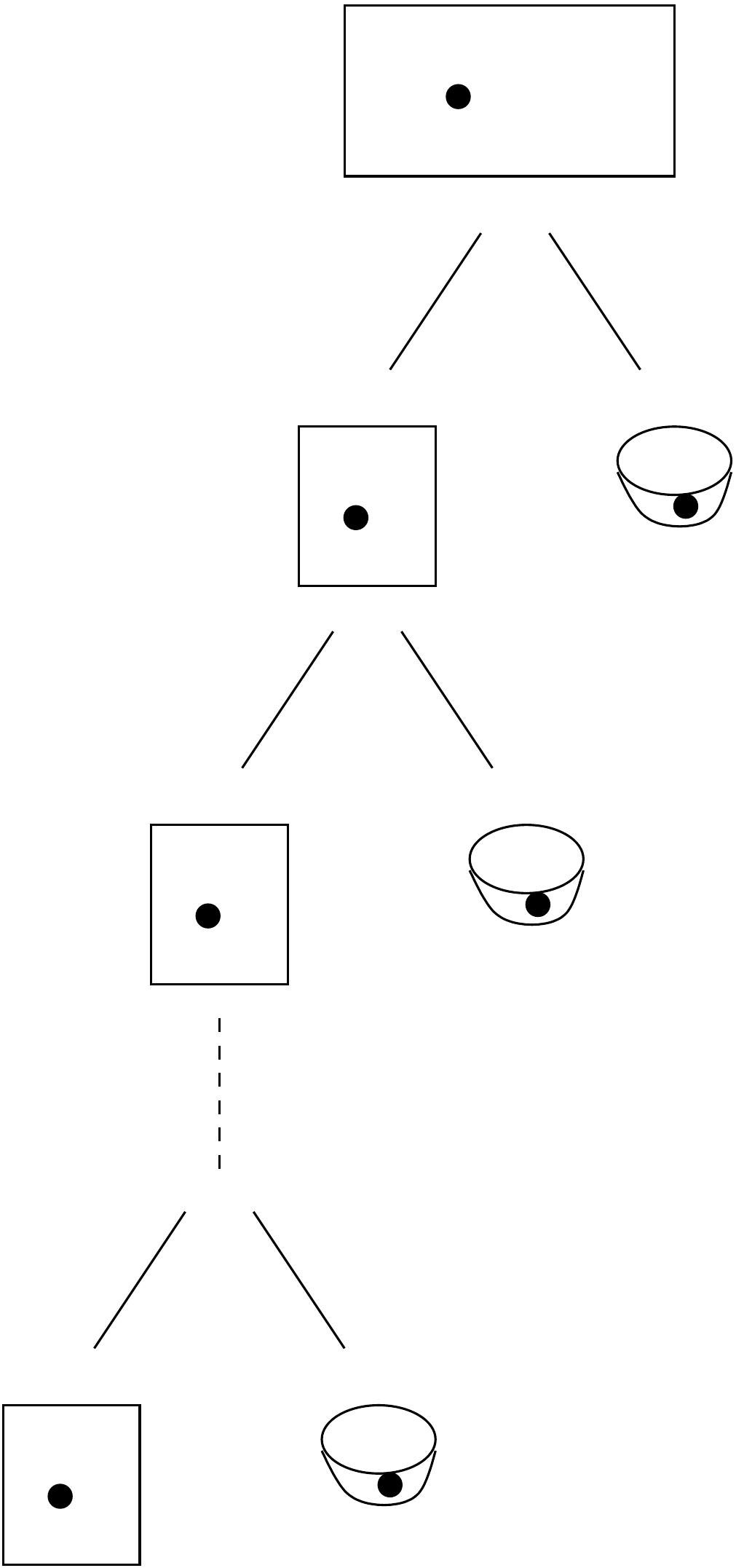}  \hspace{-3.8 in} \raisebox{.05 in}{Incompressible non-disk, non-sphere $\rightarrow$}
\end{center}

It can be seen in the figure that eventually an incompressible non-disk and non-sphere is arrived at with a dot, which is trivial in the quotient.  This can be traced all the way back up the tree to see the original surface is also trivial in the quotient.

\end{proof}

\begin{lemma}
\label{pairofpants}
Let S be a foam.  If placing a bridge on S turns two essential boundary curves into one essential boundary curve, or the placing of a bridge on S turns one essential boundary curve into two essential curves then the result of placing this bridge on S is a foam that is trivial in the quotient.  

\end{lemma}

\begin{proof}
By the relations we may assume we are starting with vertical annuli.  Thus if two essential boundary curves are turned into one essential boundary curve, we are left with a component with three essential curves and euler characteristic of negative one.  Note a surface that has three boundary components and euler characteristic negative one is a disk with two holes.  Since all of the boundary components are essential, this twice-punctured disk is incompressible.  Then since it has negative euler characteristic it is trivial in the quotient.

\end{proof}

\begin{lemma}
\label{nonhomo}

If two non-homotopic essential curves are bridged together, then the new foam is trivial in the quotient.

\end{lemma}

\begin{proof}

If two non-homotopic essential boundary curves are bridged together, then the result is one essential boundary curve since it could only be inessential if the original curves were homotopic.  Thus by Lemma \ref{pairofpants}, this new foam is trivial in the quotient.

\end{proof}

\begin{example}

It may seem that the order of applying $d_p$ and applying the neck-cutting may affect the orientation on the boundary curves.  This is an example of one way this can happen.

\begin{center}\includegraphics[height = 6 in]{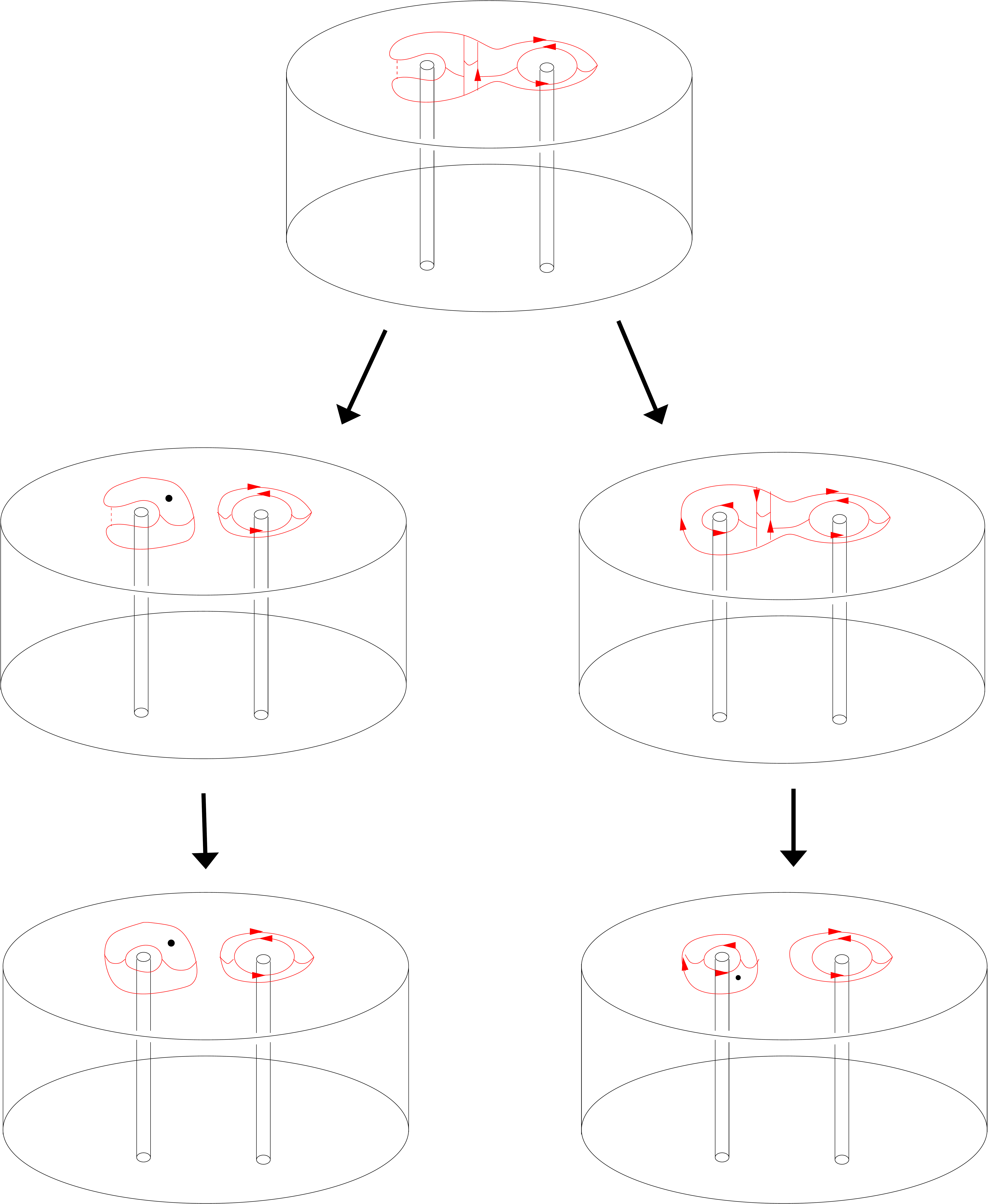}\end{center}

Start with the foam in the center.  The path to the right places a bridge first, then compresses and the path to the left compresses, then places a bridge.  Note the orientation on some curves differs on the final surface, however in each resulting foam, there is an incompressible annulus with a dot, so both are trivial, and therefore equal in the quotient.

\end{example}

The following lemma shows that this is what happens in general.

\begin{lemma}
\label{ncwelldef}

The boundary operator together with the neck-cutting relation is well-defined with regard to orientation on the quotient.

\end{lemma}
\begin{proof}

Note that on the left side of the neck-cutting relation boundary orientation from essential boundary circles can be forced on other essential boundary circles through the neck.  However, on the right side of the relation all essential boundary curves do not need to be compatible with one another since they may no longer lie on the same connected component.

\begin{center}
\includegraphics[height = .8 in]{BNrel.pdf}
\end{center}

Consider that this only becomes a problem if each component has an essential boundary component on the right-hand side.  Otherwise the compatibility of boundary curves is not an issue.  Then note one component has a dot in each summand, so by Lemma \ref{nontrivdot} both sides are trivial in the quotient.  Thus orientations may differ, but the foams are trivial, and therefore equal in the quotient.  

\end{proof}

\begin{theorem}
The boundary operator is well defined on the quotient.
\end{theorem}
\begin{proof}
Note the boundary operator is defined on the original module but it actually operates on a quotient of that module.  Thus it needs to be verified that two representations of the same class go to the same class under the boundary operator.  Thus assume $[S] = [S']$ in our quotient.  Therefore $S-S' \in B$.  So we must show that $d(S) - d(S') \in B$.  Since $d$ is linear this is equivalent to showing $d(S-S') \in B.$  Thus it is sufficient to show that given $b \in B$, that $d(b)\in B$.

Cases:
\begin{enumerate}

\item If the boundary operator is applied to surfaces that are related by the neck-cutting relation, then the result is equal in the quotient.

\begin{proof}

It must be shown that if foams are related by the neck-cutting relation before applying the boundary operator they are related after applying the boundary operator as well.  Placing a bridge does not remove any compressing disks, so the only issue is how orientations are affected.  Lemma \ref{ncwelldef} shows if the orientations agree before applying the boundary operator, they agree after applying it as well.  Thus the boundary operator is well-defined with respect to the neck-cutting relation.

\end{proof}

\item If the boundary operator is applied to a foam that has an incompressible component that isn't a disk or a sphere, but has a dot, it remains trivial in the quotient.

\begin{proof}

Let $S$ be a foam represented by a surface that has an incompressible non-disk component (thus this component has euler characteristic less than or equal to zero and has a dot).  If the boundary operator does not affect this component then it is clearly still trivial in the quotient, so we may assume we bridge this component to itself or another incompressible non-disk component.  (The disk case is immediate since if this component is bridged to a disk what remains is still an incompressible non-disk non-sphere component with a dot.)  

After placing the bridge the new connected component has euler characteristic less than or equal to negative one and it has a dot. If this new component has an essential boundary component then it is trivial in the quotient by Lemma \ref{ncwelldef}.  

If this component has no essential boundary components then it must have only one boundary component and it is inessential since we started with no inessential boundary circles and only one inessential circle can be created by a bridge.  Note there is a compressing disk present since this surface has an inessential boundary circle but is not a disk by noting the euler characteristic.  This compressing disk is separating and compressing upon it yields a disk and a closed surface with a dot.  This closed surface is not a sphere since there was a compressing disk.  Thus it a higher genus surface which has a non-separating compressing disk in an I-bundle.  Compressing upon this disk yields a surface with two dots which is trivial in the quotient.

\end{proof}

\item If the boundary operator is applied to any foam that has an incompressible component with negative euler characteristic it remains trivial in the quotient.

\begin{proof}

Let $S$ be a foam represented by a surface that has an incompressible component with negative euler characteristic.  We may assume the bridge is placed on this component.  Now assume we are bridging to another incompressible component (but not a disk) or itself.  Note before bridging these components there is at least three essential boundary components between them and after bridging there is at least one essential boundary component left.  
 
If this surface is incompressible then we are done as the euler characteristic is still negative.

If the surface is compressible compress.  If the compressing disk is non-separating then by Lemma \ref{nontrivdot}, the surface is trivial in the quotient.  If the compressing disk is separating, then after compressing we obtain two surfaces one of which has an essential boundary component.

\begin{center}\includegraphics[height = .5 in]{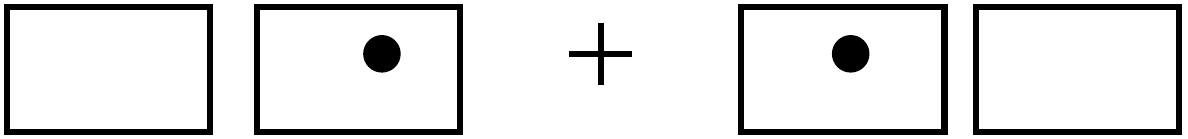}\end{center}

 \hspace{ .8 in}   essential $\partial$  \hspace{1.9 in}  essential $\partial$

The second summand drops out, by Lemma \ref{nontrivdot}, so in either case we are left with a component with an essential boundary component.  If the component with a dot has an essential boundary curve as well then the entire surface is trivial in the quotient by Lemma \ref{nontrivdot}.  

Thus we may assume the component with a dot in the first summand does not have an essential boundary curve.  If the component with a dot does not have an inessential boundary component note it is not a disk or a sphere, but it has a dot, so the entire surface is trivial in the quotient. 

 If that component does have an inessential boundary component, note at most one inessential boundary component can be created when a bridge is placed.  Thus that component has at most one inessential boundary component.  It that component is not a disk we may compress to produce a disk and a closed surface.  This closed surface has a dot, so as previously noted, it is zero in the quotient.  

Thus we may assume the second component with a dot is a disk.  If the component with essential boundary is incompressible we are done since it has negative euler characteristic.  If not, then compress.  

If the compressing disk is non-separating, then we are left with a component with an essential boundary curve and a dot, so the foam is trivial in the quotient by Lemma \ref{nontrivdot}.

If the compressing disk is separating then after we compress we have two components.  Neither are a sphere or a disk and thus they are planar surfaces with essential boundary components or closed surfaces of genus greater than or equal to one.  They are in a sum where one of them has a dot in each summand.  

\begin{center}\includegraphics[height = .5 in]{proofschm2.pdf}\end{center}

By Lemma \ref{nontrivdot} the planar surface with a dot is equal to zero.  As noted before closed non-sphere surface with dots are trivial in the quotient.  Thus after placing a bridge the surface remains trivial in the quotient.

\end{proof}

\item If the boundary operator is applied to both sides of the UTA relation then the results are equivalent in the quotient.

\begin{proof}

Based on the previous two sections of this proof, we do not need to consider if we bridge to an incompressible component with negative euler characteristic or anything that isn't a disk, but has a dot, as these cases have already been addressed.  Thus we may consider that to begin with all components are incompressible annuli without dots or disks with one or no dots.

There are four cases:

\begin{enumerate}
\item The annulus on the left side of the UTA relation is bridged to itself.

If one curve becomes two essential curves, both sides are trivial in the quotient by Lemma \ref{pairofpants}.  

If one curve becomes one essential curve, and one inessential curve, then after compressing we end up with a disk with a dot and the annuli we started with.

If two curves become one, they must become inessential since they are parallel, so then the left hand side is:

\includegraphics[height = .75 in]{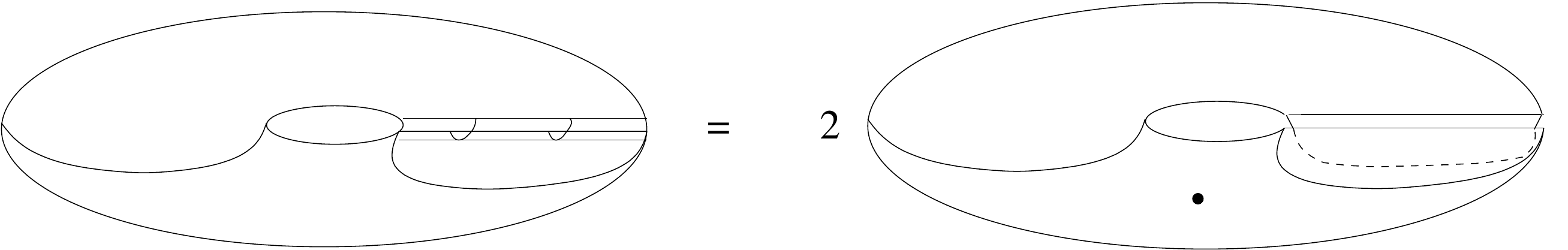}

and the right hand side is:

\includegraphics[height = .75 in]{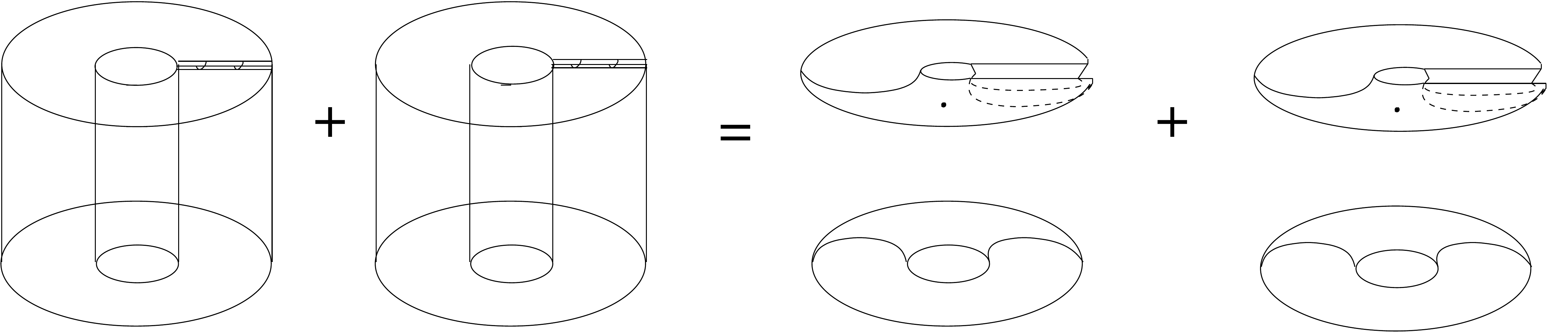}

which are equivalent in the quotient.

\item The annulus on the left side of the UTA relation is bridged to disk.

This case is immediate since bridging an annulus to a disk results in an annulus isotopic to the original annulus.

\item The annulus on the left side of the UTA relation is bridged to an unoriented top annulus.

If the annuli are not nested, we get zero on both sides by Lemma \ref{nonhomo}, so assume annuli are nested.  

The left side becomes:

\includegraphics[height = .75 in]{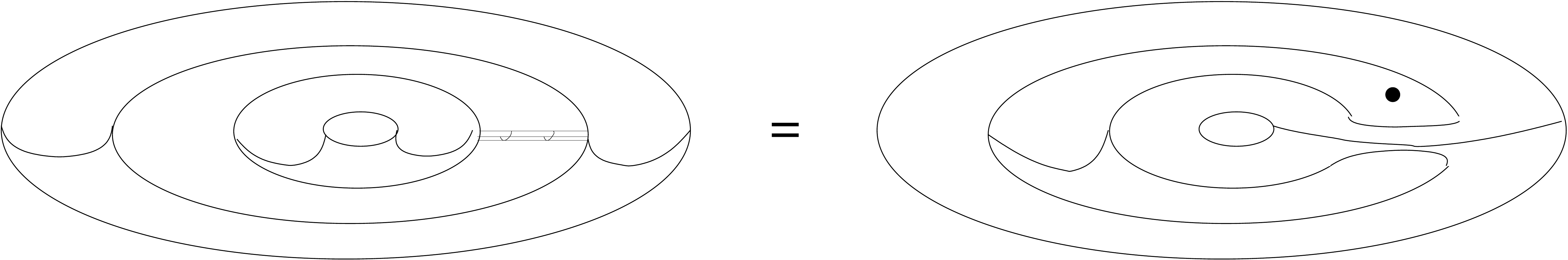}

and the right side becomes:

\includegraphics[height = .7 in]{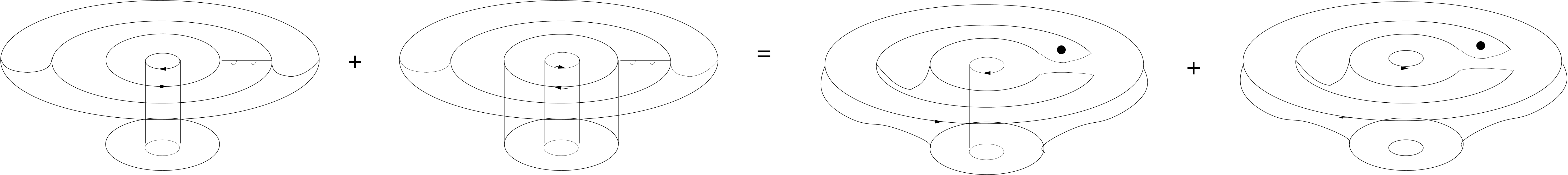}

which are equivalent in the quotient.

\item The annulus on the left side of the UTA relation is bridged to oriented vertical annulus.

If the annuli are not nested, we get zero on both sides again by Lemma \ref{nonhomo}.

Thus, assume annuli are nested.

The left side becomes:

\includegraphics[height = .75 in]{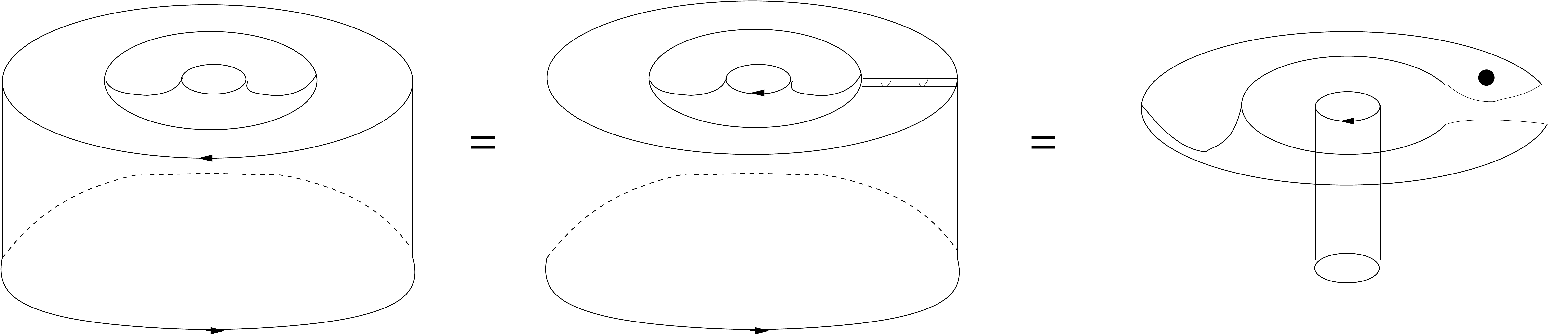}

and the right side becomes:

\includegraphics[height = .75 in]{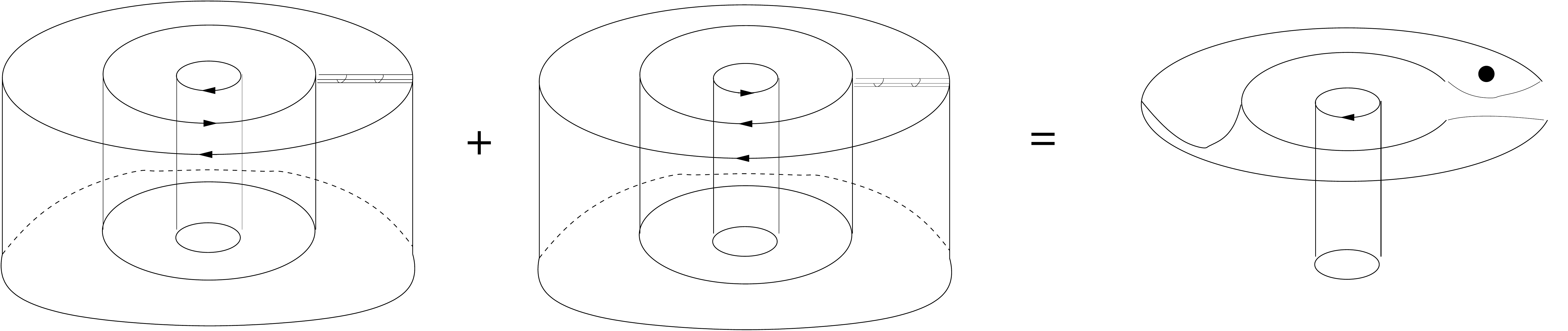}

which is are equal.

\end{enumerate}

\end{proof}

\item If the boundary operator is applied to a foam that has an annulus component with its boundary completely in the bottom and the boundary operator is applied to a foam with that annulus removed the results are equivalent in the quotient.

A bridge can't be placed in the bottom, thus the annulus will not affect the boundary operator.  Therefore the resulting foams will still be equivalent in the quotient as everything else will be affected in the same manner under the boundary operator.

\end{enumerate}

Thus when the boundary operator is applied to elements of $B$, they remain in $B$, so by the remarks at the beginning of this proof, the boundary operator is well defined on the quotient.

\end{proof}

\section{$d^2 = 0$}

\begin{lemma}

Let a and b be two crossings of a given diagram and let S be a foam.  If $d_a(d_b(S)) = 0$ because of conflicting orientations on boundary curves (EO), then $d_b(d_a(S))$ is trivial in the quotient.
\end{lemma}

\begin{proof}
Assume that any curves being combined are homotopic by Lemma \ref{nonhomo}.   Otherwise the result would be trivial in the quotient, by the previous lemma.  Since $d_a(d_b(S)) = 0$ by (EO) then there must be two components with the same orientation present. Thus we must be dealing with at least two components, and note we have at most three components since we only are dealing with two crossings.  

If there are three components then it is possible that the boundary curves are all oriented the same way or the middle component is oriented the same way as one of the other two components.  If they are all oriented the same way consider the figure:

\begin{center}
\includegraphics[height = .75 in]{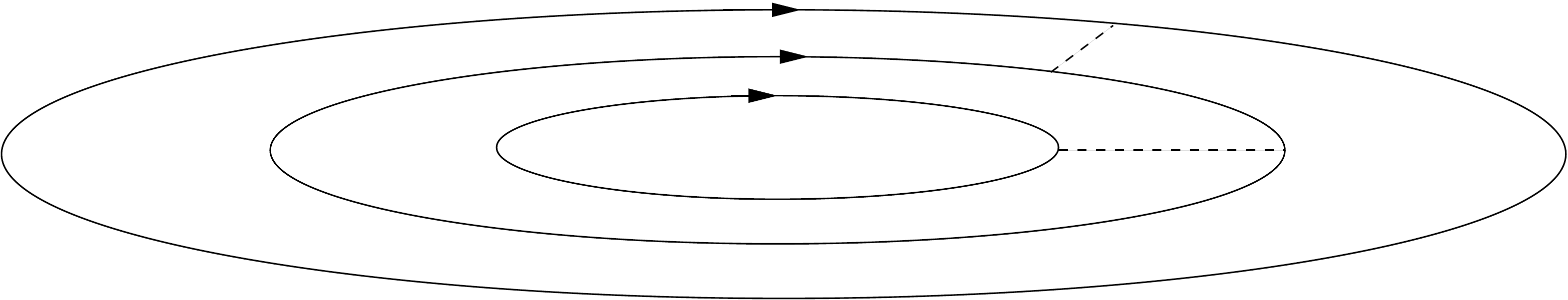}
\end{center}

Both bridges bridge two components together so the dotted lines can represent where the bridges will go.  One can see that placing either bridge results in (EO), so the order does not matter. 

Thus assume the middle component is oriented the same way as one of the other two.  Consider the following figure:

\begin{center}
\includegraphics[height = .75 in]{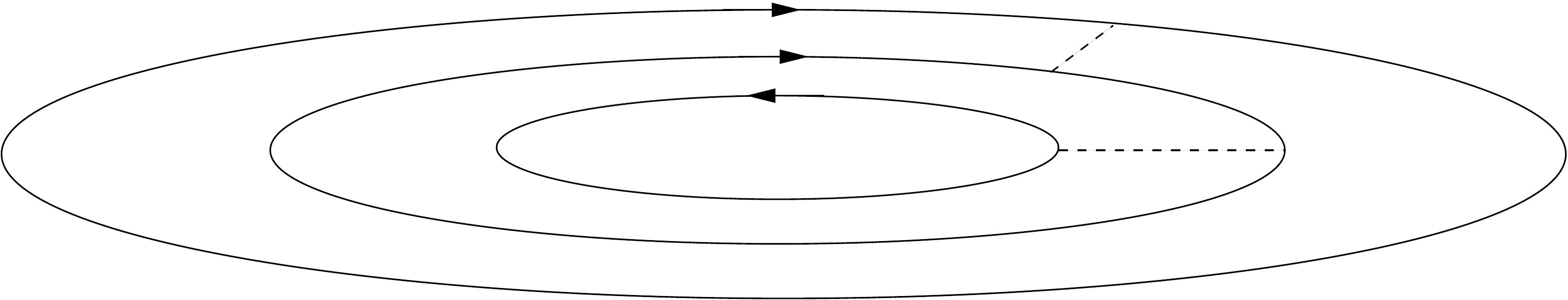}
\end{center}

One bridge may be placed, so assume we place that bridge.  That results in an annulus and a disk with a dot as in this figure:

\begin{center}
\includegraphics[height = .75 in]{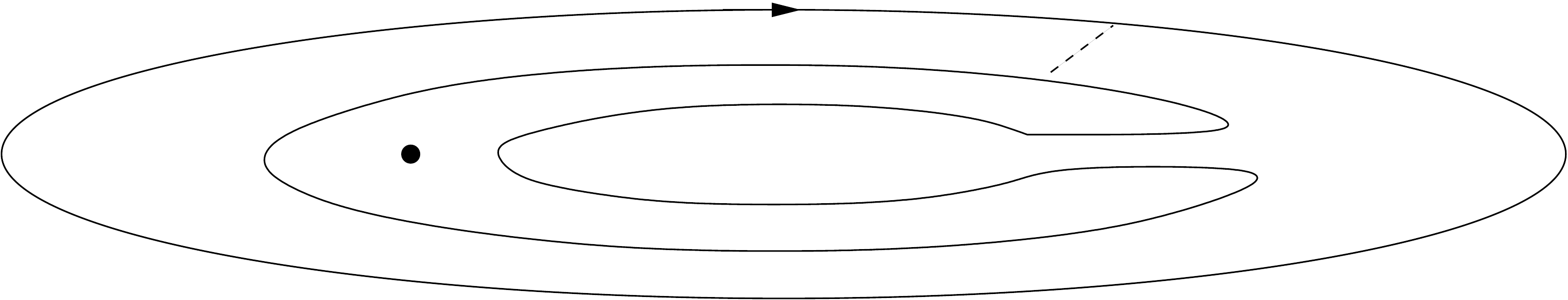}
\end{center}

This results in the second bridge connecting a non-disk and a disk with a dot, which results in a foam that is trivial in the quotient.

If there are two components then both components are oriented annuli, oriented the same way.  These annuli are connected by at least one crossing.  

\begin{center}
\includegraphics[height = .75 in]{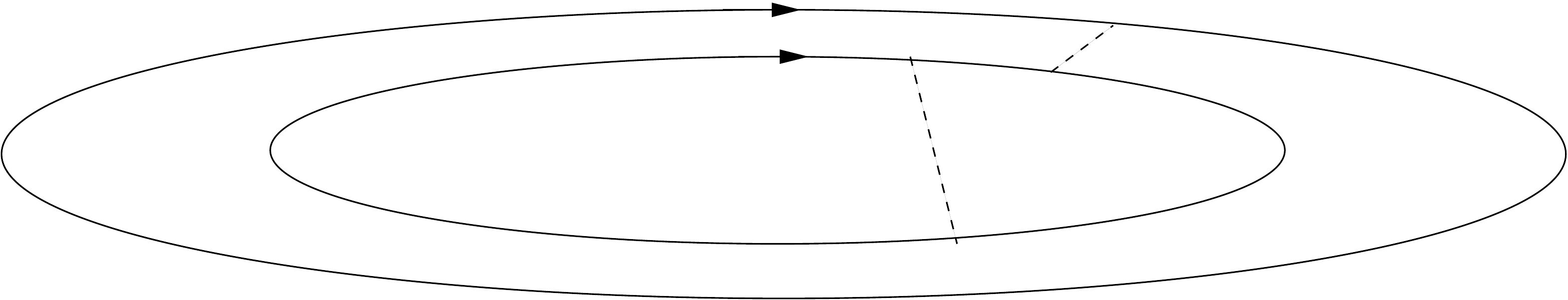}
\end{center}

If this crossing is changed first we get zero because of conflicting boundary orientation.  If another crossing is changed first the arcs either retain their original orientation, or become curves on a surface with a dot after compression.  

\begin{center}
\includegraphics[height = .75 in]{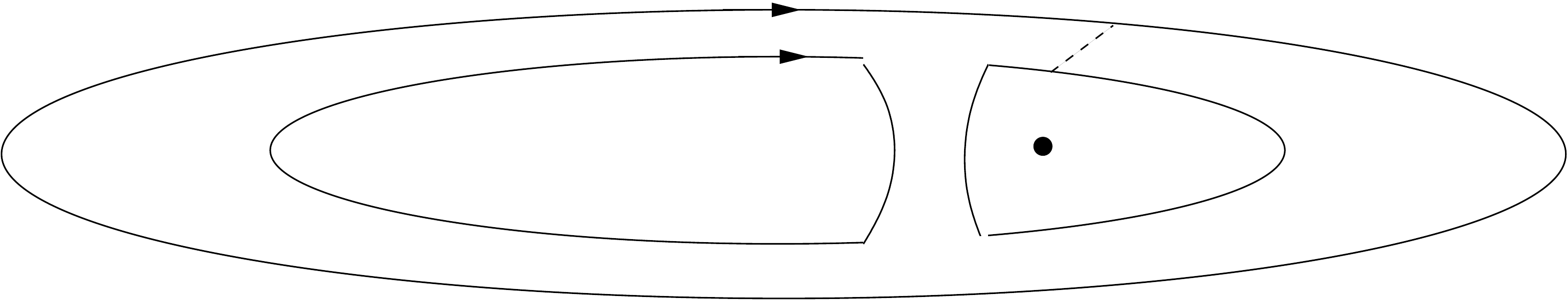}
\end{center}

Either way they can not combine with the annulus that wasn't changed which results in a foam that is trivial in the quotient.

\end{proof}

\begin{lemma}
Let a and b be two crossings of a diagram.  If a non-orientable surface is created when applying the boundary operator in the order $d_a(d_b(S))$, then $d_b(d_a(S))$ is trivial in the quotient.
\end{lemma}
\begin{proof}

Assume we start with 3 boundary circles in the top.  If we allow non-orientable surfaces to be created then placing a bridge can either make the number of boundary circles on top go up by 1, go down by 1, or stay the same.  Thus 3 circles can become 2, 3 or 4 circles after placing a bridge, which can become 1, 2, 3, 4 or 5 after placing two bridges.  

It we do not allow non-orientable surfaces to be created, then 3 circles can become 2 or 4 circles after placing a bridge, which can become 1, 3 or 5 circles after placing another bridge.

Note the only way they can agree at the end is if we end up with 1, 3 or 5 circles.

If 1 or 5 circles are left we changed the number of circles at each crossing, so no non-orientable surfaces were created while bridging.  Thus the only possibility to have created non-orientable surfaces is if we end up with 3 circles.  Note if (NOS) occurrs at all, then to get 3 circles it must have occurred at both crossings.

Also note by Lemma \ref{nonorsurf} (NOS) only occurs if the number of boundary curves stays constant.  However, if two boundary curves are bridged together the number of boundary curves always goes down.  Thus one boundary curve must become one boundary curve after bridging and this can only happen if we are bridging one curve to itself.  Thus we only need at most 2 top boundary circles to start with.  However if we have two circles and they each are connected to themselves, then clearly order of crossing change does not matter.  Therefore we can assume we start with only 1 circle and end with 1 circle.  Thus the ordering where (NOS) doesn't occur starts with 1 circle, there are 2 circles after bridging once, and 1 circle again after bridging twice.

Note if the initial boundary circle bounds a disk and placing a bridge doesn't change the number of boundary circles, then the circle bounds a non-orientable surface by the previous lemma.  Also this surface has boundary only in the top.  Since this surface was obtained from a disk by placing a bridge it has euler characteristic zero and note it has one boundary component thus it is a mobius band.  Now reflect this surface across $F$ x $\{0\}$ so that there are two mobius bands meeting only at their boundary components that lie in $F$ x $\{0\}$.  Thus we have a Klein bottle which is embedded in $\mathbb{R}^3$ which is not possible.

Thus we must conclude that the circle starts out as the boundary of a vertical annulus.  Note when the bridge is placed, the curve is split and we either get two essential curves or an inessential curve and an essential curve.    If two essential curves occur the Lemma holds by Lemma \ref{pairofpants}.  If the inessential curve and essential curve case occurs then the result after compressing is a copy of the original annulus and a disk with a dot.  

Note the annulus and the disk with a dot must recombine, but this results in an incompressible annulus with a dot, which is trivial in the quotient.

\end{proof}

%\begin{lemma}
%If placing a bridge results in a 0 based on orientation, then placing any other bridge first and then that bridge also results in a 0.  i.e. the differential commutes in these situations.

\begin{lemma}

Let $S$ be a foam and $a$ and $b$ are crossings in the associated diagram.  Then ${d}_a({d}_b(S))$ has the same orientation on boundary curves as $d_b(d_a(S))$. 

\end{lemma}
\begin{proof}
Assume $S$ has some oriented boundary components and $S$ is represented by vertical annuli and disks.  Place the bridges on $S$, but don't apply the neck-cutting relation.  Note the oriented curves on the bottom of $S$ are not affected by placing bridges and the only way $S$ can have oriented curves is if there are oriented curves in the bottom.  

After placing the bridges we are dealing with one connected component, since we haven't applied the neck-cutting relation.  After placing the bridges in both orders any oriented essential curves are necessarily oriented the same way since they must be oriented compatibly with the curves on the bottom of the original surface.

\end{proof}

\begin{theorem}
$d^2=0$

\end{theorem}

\begin{proof}
Note that by how the negative signs are distributed in the definition of the boundary operator all that needs to be shown is that given two crossings $a$ and $b$ and a foam $S$, that $d_a(d_b(S)) = d_b(d_a(S))$.  This is clearly the case if (EO) and (NOS) do not occur.

In these cases we have shown by the preceding lemmas that if they occur for one order of a and b, then  $d_a(d_b(S))$ and $d_b(d_a(S))$ are both trivial in the quotient.

Thus the partials always commute, so after adding in the appropriate negative signs, $d^2=0$.

\end{proof}

\section{Equating the Homology Theories}

Let $p$ be a crossing of the diagram $D$.  Consider the skein triple in $F$:

\medskip

\begin{center}
\begin{tabular}{ccccc}

\includegraphics[height = 1 in]{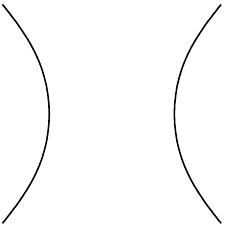}  & \hspace{.5 in} & \includegraphics[height = 1 in]{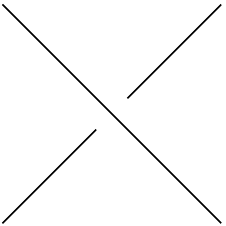}  & \hspace{.5 in} & \includegraphics[height = 1 in]{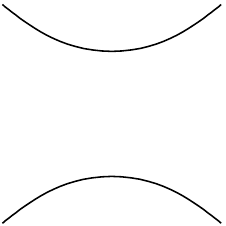} \\

        $D_{\infty}$    & &     $D_p$     & &   $D_0$   \\
 \end{tabular}

\end{center}

Now define:

$\alpha_0:C_{i,j,s}(D_{\infty}) \rightarrow C_{i-1,j-1,s}(D_p)$ is the the natural embedding as depicted in the figure below.

\begin{center}
\includegraphics[height = .75 in]{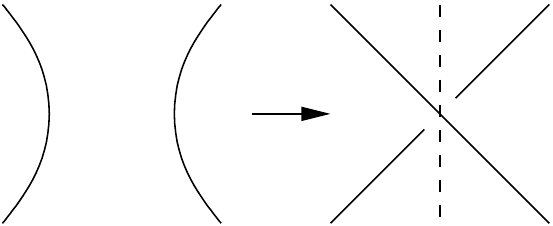}
\end{center}

$\beta:C_{i,j,s}(D_{p}) \rightarrow C_{i-1,j-1,s}(D_0)$ is the natural projection where foams with a negative smoothing at $p$ are sent to 0 and foams with a positive smoothing are affected as in the figure below.

\begin{center}
\includegraphics[height = .75 in]{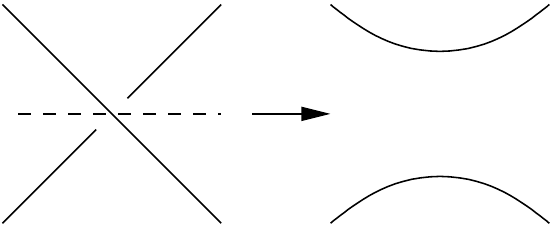}
\end{center}

Let $\alpha:C_{i,j,s}(D_{\infty}) \rightarrow C_{i-1,j-1,s}(D_p)$ be defined by $\alpha(S) = (-1)^{t'(S)}\alpha_0(S)$ where \\ $t'(S)= |$\{ $j$ a crossing of D : $j$ is before $p$ in the ordering of the crossings of $D$ and $j$ is smoothed negatively in boundary state of S\}$|$

\begin{theorem}

$\alpha$ and $\beta$ are chain maps and the sequence 

$0 \rightarrow C_{i+1,j+1,s}(D_\infty) \xrightarrow{\alpha} C_{i,j,s} (D_p)\xrightarrow{\beta} C_{i-1,j-1,s}(D_0) \rightarrow 0$

is exact.

\end{theorem}

\begin{proof} (from [APS])

First note that $d_p(\alpha(S)) = 0$ since $\alpha(S)$ has a negative smoothing at the $p$-th crossing, so the boundary operator in the $p$ direction is always 0.  

Let $\hat{d}_{q} = (-1)^{t(S,q)}d_{q}$.

Let $S_q = d_q(S)$. Note that if $q \neq p$ then $d_q(S)$ is topologically the same in  $C_{i,j,s}(D_\infty)$ and $C_{i,j,s} (D_p)$, thus $\alpha_0$ and $d_q$ commute.

Also note that either $t'(S) + 1 = t'(S_q)$ or $t'(S) = t'(S_q)$ depending on whether $p$ or $q$ comes first in the ordering of crossings.  Consider that if $t'(S) + 1 = t'(S_q)$ then $t(\alpha_0(S),q) + 1 = t(S,q)$ and if $t'(S) = t'(S_q)$, then $t(\alpha_0(S),q) = t(S,q)$.  So, in either case $t(\alpha_0(S),q) + t'(S) = t'(S_q) + t(S,q)$ mod 2.

Then note 

\begin{equation*} 
\hat{d}_{q} (\alpha(S)) = \hat{d}_q((-1)^{t'(S)}\alpha_0(S))=(-1)^{t'(S)}\hat{d}_{q}(\alpha_0(S))
\end{equation*}

\begin{equation*} 
=(-1)^{t'(S)}(-1)^{t(\alpha_0(S),q)}{d}_{q}(\alpha_0(S))=(-1)^{t'(S)+t(\alpha_0(S),q)}{d}_{q}(\alpha_0(S))
\end{equation*}

\begin{equation*} 
=(-1)^{t'(S_q)+ t(S,q)} \alpha_0({d}_{q}(S))=(-1)^{t'(S_q)} \alpha_0((-1)^{t(S,q)}{d}_{q}(S))=\alpha(\hat{d}_{q}(S))
\end{equation*}

Then we have $d(\alpha(S))= \sum_{q\neq p} \hat{d}_{q}(\alpha(S)) = \alpha(\sum_{q\neq p} \hat{d}_{q}(S))=\alpha(d(S))$, so $\alpha$ is a chain map.

\medskip

Now consider $\beta$.  Note $\beta(d_p(S)) = 0$ since $d_p(S)$ is smoothed negatively at $p$ and $\beta$ sends foams with the negative smoothing at $p$ to 0.  Also, for $q\neq p$ $d_q(\beta(S))=\beta(d_q(S))$ since the bridge is placed away from $p$, so the result is the same.  Also, $\beta$ doesn't change the number or placement of negative crossings, so we have $\hat{d}_q \beta = \beta \hat{d}_q$.  Then $d(\beta(S)) = \beta(d(S))$ and thus $\beta$ is a chain map.

\medskip

Now the exactness of the sequence is addressed.  Since $\alpha$ is an embedding it is 1-1.  The image of $\alpha$ is all foams in $C_{i,j,s} (D_p)$ that have a state as the top boundary smoothed negatively at $p$.  The kernel of $\beta$ is precisely these foams.  Since $\beta$ is a projection, it is onto $C_{i-1,j-1,s}(D_0)$.  Thus the sequence is exact.

\end{proof}

Let $\bar{C}_{i,j,s}(D)$ be the chain groups defined in [APS] and $\bar{H}_{i,j,s}(D)$ be the homology groups defined in [APS].   Let $H_{i,j,s}$ be the homology of the chain complex we have constructed in the previous sections.

In [APS] a circle in $F$ is said to be trivial if it bounds a disk in $F$ and non-trivial otherwise.  Thus when referring to circles of an enhanced state coming from [APS] these terms will be used.  Also note that trivial circles correspond to inessential boundary curves in the top and non-trivial circles correspond to essential circles in the top.

\begin{definition} $\Phi: \bar{C}_{i,j,s} \rightarrow  C_{i,j,s}$ is defined by taking an enhanced state in $\bar{C}$ and changing each circle as follows to get a foam with the same state in $C$:

Trivial circle marked with a + $\rightarrow$ disk with a dot

Trivial circle marked with a $-$ $\rightarrow$ disk without a dot

Nontrivial circle marked with a $+0$ $\rightarrow$ vertical annulus with the positive orientation

Nontrivial circle marked with a $-0$ $\rightarrow$ vertical annulus with the negative orientation

\end{definition}

\begin{example}
This is an example of how the $\Phi$ map affects an enhanced state from [APS].

\begin{center}
\includegraphics[height = 1 in]{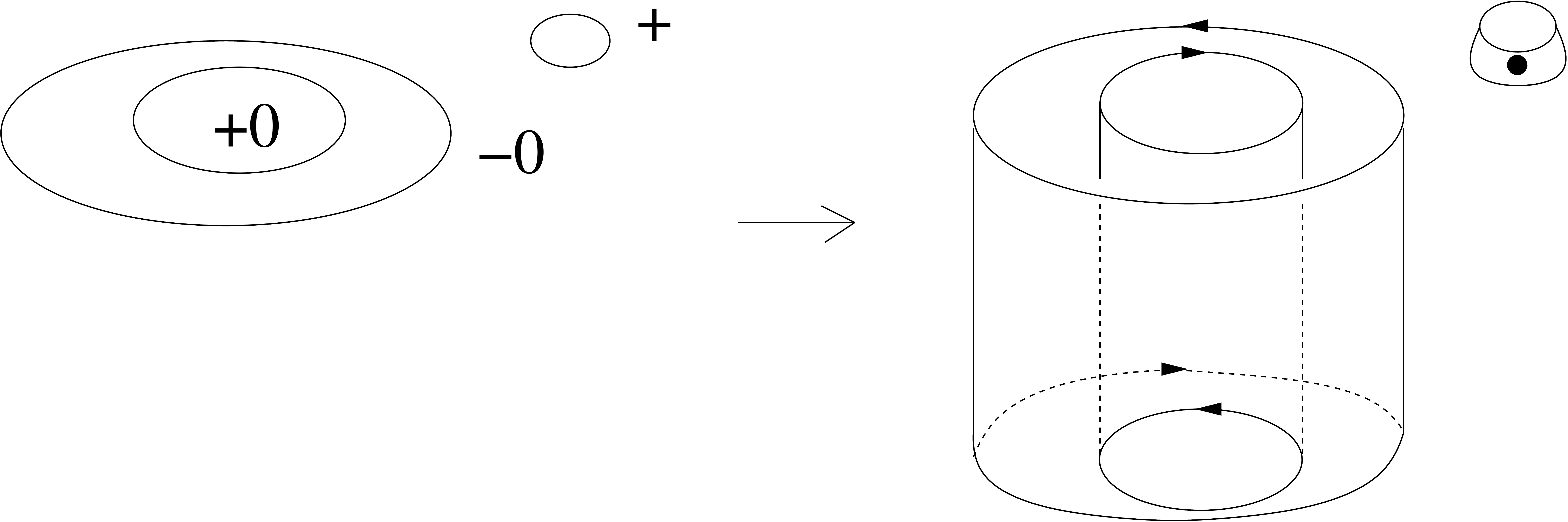}
\end{center}

\end{example}

\begin{lemma}

$\Phi: \bar{C}_{i,j,s} \rightarrow  C_{i,j,s}$ is a chain map $\forall i,j,s$.

\end{lemma}

\begin{proof}  To show $\Phi$ is a chain map, we need to show $\Phi \tilde{d}=d \Phi$. Where $\tilde{d}$ is the boundary operator coming from [APS].  First it will be shown that $\Phi \tilde{d}_i = d_i \Phi$ $\forall i$.

This is treated by cases. In the [APS] theory a trivial circle may have a + or a $-$ and a non-trivial circle may have a +0 or a $-$0.  Thus a + will refer to a trivial circle marked with a + and a +0 will refer to a non-trivial circle marked with a +0.  The notation is similar for $-$ and $-0$.  If the marking on a circle is not specified then T refers to a trivial circle and N refers to a non-trival circle.  Based on how trivial and non-trivial curves may change when a smoothing changes here are the possible cases:

\begin{center}

\begin{tabular}{|c|c|c|c|}
\hline
How curves  & Possible initial markings& Possible& Number of   \\

may change &&outcomes& possibilities\\

\hline

T $\rightarrow$ TT or NN   &   2 choices, + or $-$ &  2 results   &   $ 2*2 = 4$       \\  \hline

N $\rightarrow$ NT or NN    &  2 choices, +0 or $-0$  & 2 results  & $ 2*2 = 4$        \\  \hline

TT $\rightarrow$ T        &         $2 + 1 =3$ choices, ++, +$-$, or $--$  & 1 result  &   3           \\   \hline

 NN $\rightarrow$ T or N  &    $2 + 1 = 3$ choices, +0+0, +0$-0$, or $-0-0$  & 2 results  &   $ 3*2 = 6$   \\  \hline

 TN $\rightarrow$  N     &       $2*2 = 4$ choices, ++0, +$-0$, $-+0$ or  $--0$  & 1 result   &   $ 4$   \\   \hline
 
 \hline
 \end{tabular}

 \end{center}

 Total $= 4 + 4 + 3 + 6 + 4 = 21$ cases
 
\medskip

The boundary operator from [APS] is determined by how circles change (with respect to being trivial and nontrivial) when a smoothing is switched.  The table below shows what the partial boundary operator for the theory coming from [APS] does in all of the above cases when one crossing is switched and then $\Phi$ is applied.

\begin{center}

\begin{tabular}{|cc|l|c|c|c|c|c|}
	\hline
	& How the circles change &$\Phi \tilde{d}_i$ \\ 
	\hline
1. &  T $\rightarrow$ TT   &$\Phi(\tilde{d}_i$(+)) = $\Phi$(++) = \includegraphics[height = .3 in]{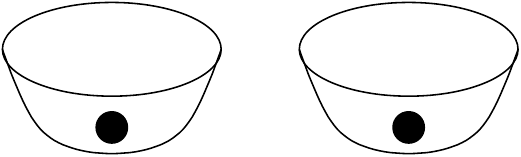}  \\
  \hline
2. &  T $\rightarrow$ NN & $\Phi(\tilde{d}_i(+)) = \Phi(0) = 0$ \\
  \hline
3. &  T $\rightarrow$ TT & $\Phi(\tilde{d}_i(-)) = \Phi((+-) + (-+)) =$ \includegraphics[height = .3 in]{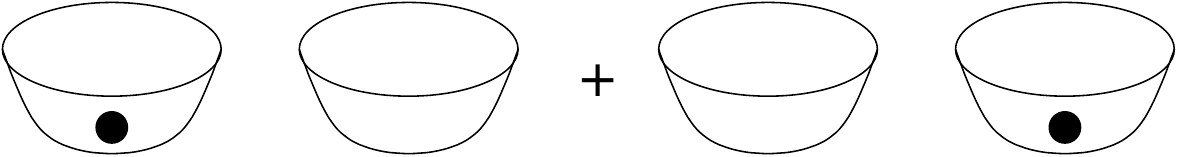} \\
\hline
4. &  T $\rightarrow$ NN &  $\Phi(\tilde{d}_i(-)) = \Phi((+0 -0) + (-0+0)) =$ \includegraphics[height = .3 in]{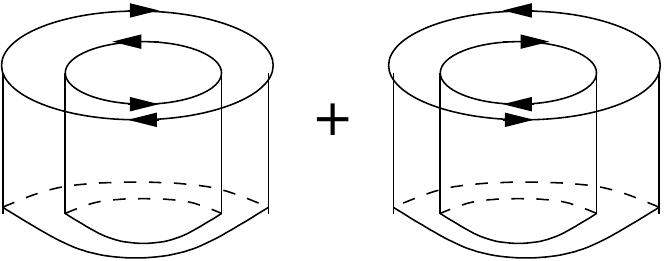} \\
\hline
5. &     N $\rightarrow$ NT &  $\Phi(\tilde{d}_i(+0)) = \Phi(+0 +) =$ \includegraphics[height = .3 in]{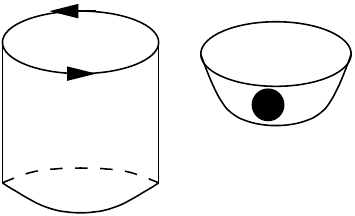} \\
\hline
6. &   N $\rightarrow$ NN & $\Phi(\tilde{d}_i(+0)) = \Phi(0) =0$  \\
\hline
7. &     N $\rightarrow$ NT &  $\Phi(\tilde{d}_i(-0)) = \Phi(-0 +) =$ \includegraphics[height = .3 in]{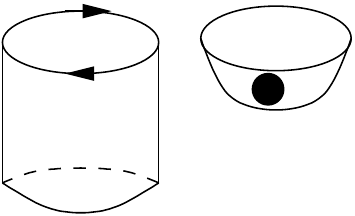} \\
\hline
8. &  N $\rightarrow$ NN & $\Phi(\tilde{d}_i(-0)) = \Phi(0) =0$  \\
\hline
9. &   TT $\rightarrow$ T  &  $\Phi(\tilde{d}_i(++)) = \Phi(0) =0$ \\
\hline
10. &   TT $\rightarrow$ T &  $\Phi(\tilde{d}_i(+-)) = \Phi(+) =$ \includegraphics[height = .3 in]{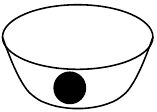} \\
\hline
11. &  TT $\rightarrow$ T &  $\Phi(\tilde{d}_i(--)) = \Phi(-) =$ \includegraphics[height = .3 in]{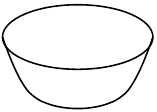} \\
\hline
12. &  NN $\rightarrow$ T &  $\Phi(\tilde{d}_i(+0+0)) = \Phi(0) =0$ \\
\hline
13. &  NN $\rightarrow$ N &  $\Phi(\tilde{d}_i(+0+0)) = \Phi(0) =0$\\
\hline
14. &   NN  $\rightarrow$ T &  $\Phi(\tilde{d}_i(+0-0)) = \Phi(+) =$ \includegraphics[height = .3 in]{cm10.pdf} \\
\hline
15. &  NN$\rightarrow$ N  &  $\Phi(\tilde{d}_i(+0-0)) = \Phi(0) =0$ \\
\hline
16. &  NN $\rightarrow$ T &  $\Phi(\tilde{d}_i(-0-0)) = \Phi(0) =0$ \\
\hline
17. &  NN $\rightarrow$ N &  $\Phi(\tilde{d}_i(-0-0)) = \Phi(0) =0$ \\
\hline
18. &  TN $\rightarrow$ N &  $\Phi(\tilde{d}_i(++0)) = \Phi(0) =0$ \\
\hline
19. &  TN $\rightarrow$ N &  $\Phi(\tilde{d}_i(+-0)) = \Phi(0) =0$ \\
\hline
20. &   TN $\rightarrow$ N &  $\Phi(\tilde{d}_i(--0)) = \Phi(-0) =$ \includegraphics[height = .3 in]{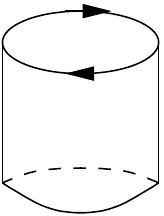} \\
\hline
21. &  TN $\rightarrow$ N &  $\Phi(\tilde{d}_i(-+0)) = \Phi(+0) =$ \includegraphics[height = .3 in]{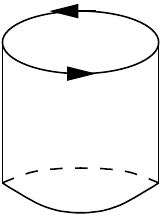} \\
\hline

\end{tabular}
	
\end{center}

Note that under $\Phi$ the associated state isn't affected, thus for example if T $\rightarrow$ TT by changing a smoothing before applying $\Phi$, then after applying $\Phi$ the boundary circles behave the same way, and an inessential boundary circle turns into two inessential boundary circles by placing a bridge.

The following 21 items show what $d_i \Phi$ is in each of the cases when the boundary circles are affected as in the previous table.

\begin{enumerate}	

\item Note $\Phi(+)$ = \includegraphics[height = .2 in]{cm10.pdf}.  After a bridge is placed there are two trivial boundary curves in the top.  This has euler characteristic equal to 0, and thus it is a compressible annulus.  Compress the annulus to get two disks, each with a dot.   %1

\item $\Phi(+)$ = \includegraphics[height = .2 in]{cm10.pdf}.  When a bridge is placed there are two non-trivial boundary components in the top.  This is an incompressible annulus with a dot, so it is trivial in the quotient. %2

\item  $\Phi(-)$ = \includegraphics[height = .2 in]{cm11.pdf}.  After a bridge is placed there are two trivial boundary curves in the top.  This is a compressible annulus.  Compress the annulus to get disk with dot, disk + disk, disk with dot.  %3

\item  $\Phi(-)$ = \includegraphics[height = .2 in]{cm11.pdf}.  After a bridge is placed there are two non-trivial boundary curves in the top.  This is an incompressible annulus, so have unoriented annulus = average of oriented annuli.  %4

\item  $\Phi(+0)$ = \includegraphics[height = .2 in]{cm21.pdf}.  After a bridge is placed there is a non-trivial boundary curve in the top and a trivial boundary curve in the top.  Compress the neck that is near the trivial boundary curve to get an annulus, oriented same way as the original annulus and a disk with a dot.  %5

\item  $\Phi(+0)$ = \includegraphics[height = .2 in]{cm21.pdf}.  After a bridge is placed there are two non-trivial boundary curves on the top.  One can only compress and separate boundary curves if we have at least 4 non-trivial and we only have three, so we have a surface that is trivial in the quotient by Lemma \ref{pairofpants} %6

\item   Refer to 5  %7

\item   Refer to 6  %8

\item    $\Phi(++)$ = \includegraphics[height = .2 in]{cm1.pdf}.  After a bridge is placed there is one trivial boundary component.  Now we have two dots on the same component, so it is trivial in the quotient.  %9

\item     $\Phi(+-)$ = \includegraphics[height = .2 in]{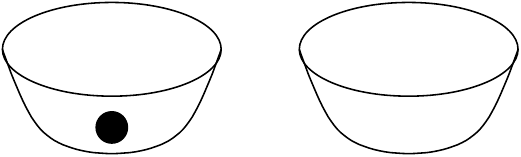}.  After a bridge is placed there is one trivial component.  These two disks combined to make a disk with a dot. %10

\item     $\Phi(--)$ = \includegraphics[height = .2 in]{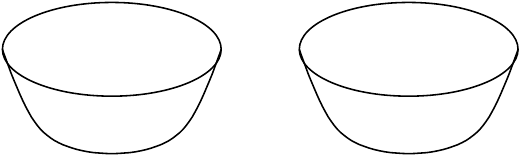}.  After a bridge is placed there is one trivial boundary component.   This leaves us with a disk. %11

\item     $\Phi(+0 +0)$ = \includegraphics[height = .2 in]{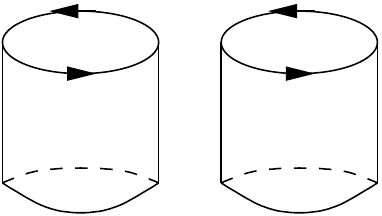}.  Placing a bridge would result in a trivial boundary component in the top.  Thus the original boundary components must have been parallel.  Therefore the bridge falls into the category of (EO) since they are oriented the same way. Thus the result is trivial in the quotient.  %12

\item     $\Phi(+0 +0)$ = \includegraphics[height = .2 in]{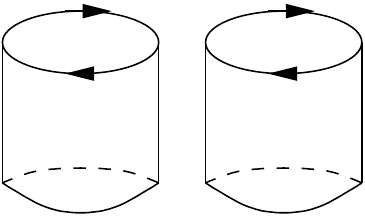}. Placing a bridge results in one non-trivial boundary curve on the top.  Thus we have an incompressible pair of pants which is trivial in the quotient.%13

\item     $\Phi(+0 -0)$ = \includegraphics[height = .2 in]{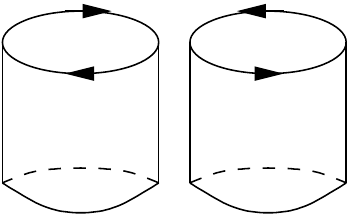}.  After placing a bridge there is one trivial boundary component. Thus the original non-trivial curves were homotopic.  Compress upon the disk that is present near the trivial curve on top.  This results in a disk on top with a dot and an annulus on the bottom + disk on top with an annulus with a dot on the bottom which is equivalent to just having a disk with a dot in the quotient. %14

\item      $\Phi(+0 -0)$ = \includegraphics[height = .2 in]{cm15dom.pdf}.  After a bridge is placed there is one non-trivial boundary component.  As in 13, we have an incompressible pair of pants which is trivial in the quotient.%15

\item   Refer to 12  %16

\item  Refer to 13  %17

\item    $\Phi(+ +0)$ = \includegraphics[height = .2 in]{cm5.pdf}.  After a bridge is placed there is one non-trivial boundary curve on the top.  Note bridging to a disk doesn't change the annulus, except it adds a dot, which makes the foam trivial in the quotient. %18

\item  Refer to 18  %19

\item    $\Phi(- -0)$ = \includegraphics[height = .2 in]{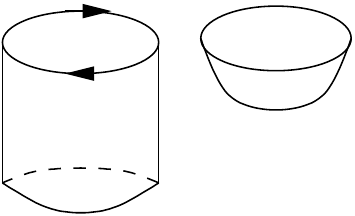}.  After a bridge is placed there is a non-trivial boundary component on top.  Absorbing a disk doesn't change annulus, so we get the same annulus with the same orientation.  %20

\item   Refer to 20  %21

\end{enumerate}

By examining the list and the table, we can see that $d_i \Phi = \Phi \tilde{d}_i$ in each case.

Thus note $\Phi (\tilde{d}(S)) = \Phi( \sum_i (-1)^{t'(S,i)}\tilde{d}_i(S) )= \sum_i (-1)^{t'(S,i)}\Phi(\tilde{d}_i(S)) = \sum_i (-1)^{t'(S,i)}d_i(\Phi(S)) = d(\Phi(S)).$  Thus $\Phi$ is a chain map, as desired.

\end{proof}

\begin{theorem}
Given a link diagram $D$, $\bar{H}_{i,j,s}(D) \cong H_{i,j,s}(D)$  $\forall i,j,s$ by $\Phi_*$.

\end{theorem}

\begin{proof}

Let $\bar{I}$, $\bar{J}$, $\bar{K}$ be the indices coming from the [APS] theory and let $\bar{S}$ be an enhanced Kauffman state from [APS].

Then note clearly $\bar{I}(\bar{S}) = I(\Phi(\bar{S}))$ since the smoothings stay the same under $\Phi$.

We also have,

$\bar{J}$ ($+$) $= \bar{I}(+)$ + 2({\# positive trivial circles} $-$ {\# negative trivial circles}) $ = I(\Phi(+)) + 2(1 - 0) = I(S) + 2(2 - 1)=I(\text{a disk with a dot}) + 2(2d - \chi(\text{a disk})) = J(\text{a disk with a dot})$

Similarly,

$\bar{J}$ ($-$) $= \bar{I}(-)$ + 2({\# positive trivial circles} $-$ {\# negative trivial circles}) $ = I(\Phi(-)) + 2(0 - 1) =I(\text{a disk without a dot}) + 2(2d - \chi(\text{a disk})) = J(\text{a disk without a dot})$

Also,

$\bar{J}$ ($+$0) $= \bar{I}(+0)$ + 2({\# positive trivial circles} $-$ {\# negative trivial circles}) $ = I(\Phi(+0)) + 2(0 - 0) = I(\text{an annulus}) =I(\Phi(+0)) + 2(2d - \chi(\text{an annulus})) = J($an annulus with bottom boundary curve oriented in the positive direction)

Finally,

$\bar{J}$ ($-$0) $= \bar{I}(-0)$ + 2({\# positive trivial circles} $-$ {\# negative trivial circles}) $ = I(\Phi(-0)) + 2(0 - 0) = I(\text{an annulus}) =I(\Phi(-0)) + 2(2d - \chi(\text{an annulus})) = J($an annulus with bottom boundary curve oriented in the negative direction)

Then note all non-trivial circles that are present in a smoothing of the diagram appear in the bottom of the foam, so the $\bar{K}$-grading is also preserved under $\Phi$.

\medskip

The proof will proceed by induction on the number of crossings in the diagram.

Assume $D$ has zero crossings.  Therefore the boundary maps are all the zero map.  Thus the chain groups are also the homology groups.  Note $\Phi$ is an isomorphism on the chain groups since it takes generators to generators, so it is also an isomorphism on homology in this case.

Let $D$ be a diagram in $F$, with $n$ crossings and inductively assume $\Phi_*$ is an isomorphism for all diagrams with less than $n$ crossings.  

Note we have a relation between the short exact sequences coming from the two theories.  The diagram commutes since $\alpha$ and $\bar{\alpha}$ are defined identically and the same is true for $\beta$ and $\bar{\beta}$.

\begin{diagram}
0 & \rTo & \bar{C}_{i+1,j+1,s}(D_\infty) & \rTo^{\bar{\alpha}} & \bar{C}_{i,j,s}(D_p)  &  \rTo^{\bar{\beta}} &  \bar{C}_{i-1,j-1,s}(D_0)  &  \rTo &  0  \\
 & & \dTo_{\Phi} &  \circlearrowright &  \dTo_{\Phi}  & \circlearrowright  &  \dTo_{\Phi}   && \\
0 & \rTo & C_{i+1,j+1,s}(D_\infty) & \rTo^{\alpha} & C_{i,j,s}(D_p)  &  \rTo^{\beta} &  C_{i-1,j-1,s}(D_0)  &  \rTo &  0\\
\end{diagram}

This induces the long exact sequence:

\begin{diagram}
\dots &  \bar{H}_{i+1,j-1,s}(D_0) & \rTo^{\partial}  & \bar{H}_{i+1,j+1,s}(D_\infty) & \rTo^{\bar{\alpha}_*} & \bar{H}_{i,j,s}(D_p)  &  \rTo^{\bar{\beta}_*} &  \bar{H}_{i-1,j-1,s}(D_0)  &  \rTo^{\partial} &  \bar{H}_{i-1,j+1,s}(D_\infty) & \rTo & \dots \\
& \dTo{\Phi_*}&  \circlearrowright & \dTo_{\Phi_*} &  \circlearrowright &  \dTo_{\Phi_*}  & \circlearrowright  &  \dTo_{\Phi_*}   & \circlearrowright &  \dTo_{\Phi_*}  \\
\dots &  {H}_{i+1,j-1,s}(D_0) &\rTo^{\partial}  & H_{i+1,j+1,s}(D_\infty) & \rTo^{\alpha_*} & H_{i,j,s}(D_p)  &  \rTo^{\beta_*} &  H_{i-1,j-1,s}(D_0)  &  \rTo^{\partial} &  H_{i-1,j+1,s}(D_\infty)  & \rTo & \dots\\
\end{diagram}

All $\Phi_*$, except the middle one, are isomorphisms by the inductive assumption.  Also, the diagram commutes since $\Phi$ is a chain map.

\medskip

Note by the five lemma the middle $\Phi$ is an isomorphism.  Thus by induction given a link diagram $D$, $\bar{H}_{i,j,s}(D) \cong H_{i,j,s}(D)$  $\forall i,j,s$ by $\Phi_*$.

\end{proof}

 Since Asaeda, Przytycki and Sikora proved invariance for the $\bar{H}(D)$ homology and $\bar{H}(D)\cong H(D)$ by the previous theorem we obtain,
 
\begin{corollary}

$H(D)$ is an invariant under Reidemeister moves 2 and 3 and a Reidemeister 1 move shifts the indices in a predictable way.

\end{corollary}

\section{References}

\begin{itemize}

\item[APS]  Asaeda, Marta M., Jozef H. Przytycki and Adam S. Sikora,  \textbf{Categorification of the Kauffman bracket skein module of I-bundles over surfaces},  Algebr. Geom. Topol. 4 (2004) 1177-1210

\item[BN1] Bar-Natan, Dror, \textbf{Khovanov's homology for tangles and cobordisms}, Geom. Topol. 9(2005) 1443-1499

\item[BN2] Bar-Natan, Dror, \textbf{On Khovanov's categorification of the Jones polynomial}, Algebr. Geom. Topol. 2 (2002) 337-370  MR1917056

\item[K] Khovanov, M., \textbf{A categorification of the Jones polynomial}, Duke Math. J. 101 (2000) 359-426 MR1740682

\item[M] Manturov, V., \textbf{Additional Gradings in Khovanov Homology}, arXiv:0710.3741v1 [math.GT]

\item[T] Turaev, Vladimir and Paul Turner, \textbf{Unoriented topological quantum field theory and link homology}, arXiv:math/0506229v2

\item[V] Viro, Oleg, \textbf{Remarks on definition of Khovanov homology}, arXiv:math/0202199v1 [math.GT]

\end{itemize}

\end{document}